\documentclass[11pt]{article}

\usepackage{graphicx,latexsym,amsfonts,graphicx}
\graphicspath{{../}}

\usepackage{algorithmicx}
\usepackage{algorithm}
\usepackage{algpseudocode}

\usepackage[usenames]{color}
\usepackage{calc}
\usepackage{verbatim}
\usepackage[mathcal]{euscript} 
\usepackage{subcaption}
\usepackage{multicol}
\usepackage{geometry}
\usepackage{setspace}
\usepackage{comment}
\usepackage{amssymb}
\usepackage{mathtools} %
\usepackage{hyperref}
\usepackage{url}
\usepackage{amsthm}

\usepackage{booktabs}  %
\usepackage{multirow}  %

\newcommand{\R}{\ensuremath{\mathbb{R}}}

\newcommand{\be}{\begin{eqnarray}}
\newcommand{\ee}{\end{eqnarray}}

\newcommand{\EE}{\mathbb{E}}

\newcommand{\indep}{\perp \!\!\! \perp}
\newcommand*{\defeq}{\stackrel{\text{def}}{=}}

\newcommand{\esteban}[1]{{#1}}
\newcommand{\andrew}[1]{{#1}}

\newtheorem{theorem}{Theorem}

\newtheorem{lemma}[theorem]{Lemma}

\title{The Monge optimal transport barycenter problem}

\author{Andrew D. Lipnick${^{1}}$, Esteban G. Tabak${^{2}}$, Giulio Trigila${^{3}}$, Yating Wang${^{2}}$, Xuancheng Ye${^{2}}$, Wenjun Zhao${^{4}}$\\
\small{${^1}$ Department of Mathematics, Lafayette College, Easton, PA 18042, USA} \\
\small{${^2}$ Department of Mathematics, Courant Institute of Mathematical Sciences, New York University, New York, NY 10012, USA} \\
\small{${^3}$ Department of Mathematics, Baruch College, City University of New York, New York, NY 10010, USA} \\
\small{${^4}$ Department of Mathematics, Wake Forest University, Winston-Salem, NC 27106, USA} \\
}

\begin{document}
 \doublespacing
\maketitle

\begin{abstract}
A novel methodology is developed for the solution of the data-driven Monge optimal transport barycenter problem, where the pushforward condition is formulated in terms of the statistical independence between two sets of random variables: the factors $z$ and a transformed outcome $y$. Relaxing independence to the uncorrelation between all functions of $z$ and $y$ within suitable inner-product spaces leads to an adversarial formulation, for which the adversarial strategy can be found in closed form through the first principal components of a matrix, and the resulting pure minimization problem can be solved efficiently through flows in phase space. The methodology extends beyond scenarios where only discrete factors affect the outcome, to multivariate sets of both discrete and continuous factors, for which the corresponding barycenter problems have infinitely many marginals.  Corollaries include a new framework for the solution of the Monge optimal transport problem, a procedure for the data-based simulation and estimation of conditional probability densities, and a nonparametric methodology for Bayesian inference. 

\end{abstract}

\section{Introduction}

A central problem in the analysis of data is to estimate how a set of variables $x \in \mathcal{X}$, the \emph{outcome}, depends on a set of  covariates $z \in \mathcal{Z}$, the \emph{factors}, a dependence that can be fully characterized by the conditional distribution $\rho(x|z)$. One seeks to extract from $n$ observed data pairs $\{x_i, z_i\}$, either an evaluation procedure for $\rho$ itself or a procedure to draw samples $\{x_j^*\}$ from $\rho(x|z_*)$ for any target value $z_*$. This is particularly challenging when $z$ includes continuous components, since any particular value $z_*$ has small probability of having appeared among the $\{z_i\}$, even less of having shown up in enough observational pairs to warrant a statistical analysis based on those pairs alone.

This article develops a data-driven methodology for the estimation and simulation of conditional distributions based on the [Monge] optimal transport barycenter problem (OTBP), seeking a map
$$ y = T(x, z) \in \mathcal{X} $$
that removes from $x$ the variability that $z$ can explain, i.e. such that the random variables $y$ and $z$ are independent. 
In order not to remove any additional variability from $x$, one selects the map $T$ that minimizes the expected value of a total transportation cost
$ C = \EE_{\pi}[c(x, T(x, z))] $,
where $\pi$ is the joint distribution of $x$ and $z$. The pairwise cost function $c(x, y)$ quantifies the deformation of the data incurred by moving $x$ to $y$.
This results in an OTBP of the form
\begin{equation}\label{Barycenter}
 \min_{y = T(x, z)} \EE_{\pi}\left[c\left(x, y\right)\right] \quad \hbox{s.t.} \quad y \indep z,
\end{equation}
where the symbol $\indep$ stands for independence. We can use the solution to this problem to simulate the conditional distribution $\rho(x|z)$ for a target value $z = z_*$, extracting $n$ samples $\{x^*_i\} \sim \rho(x|z_*)$ through
$$ x^*_i = T^{-1}\left(y_i, z_*\right), \quad  y_i = T\left(x_i, z_i\right), $$
where $T$ and its inverse $T^{-1}$ are regarded as maps between $x$ and $y$ parameterized by $z$. This procedure
first removes from $x_i$ the variability attributable to $z=z_i$, and then restores that variability with $z = z_*$, while the variable $y$ stores the variability in $x$ that $z$ does not explain.

Other uses of the OTBP include the following: 
\begin{enumerate}

 \item In order to eliminate the effect of confounding variables $z$ from the data $x$, we simply move the $\{x_i\}$ to their counterpart in the barycenter, $\{y_i = T(x_i, z_i)\}$. Examples include the removal of batch effects, the consolidation of different data bases, where $z$ represents the data source and, more generally, the removal of the confounding effects of any set of variables $z$ that are not considered in the study under way. 

\item The explanatory power of the covariates $z$ can be quantified by the total cost $C$. This ranges from the extreme scenario where $z$ has no explanatory value, so $x$ is already independent of $z$, $y = x$ and $C = 0$, to the opposite extreme where all variability in $x$ can be explained by $z$, so the barycenter reduces to a single point $\bar{y}$, which maximizes $C$.
 Quantifying through $C$ the explanatory power of $z$ gives rise to a rich methodology for factor selection and discovery \cite{TT2, yang2020conditional}.
 
 \item The barycenter problem permits not only simulating but also estimating conditional densities (see section \ref{sec:CDE}), and therefore yields a model-free, non-parametric data-based procedure for Bayesian inference: given a prior distribution $\gamma_{pr}(z)$, a set of sample pairs $\{x_i, z_i\}$ drawn from an unknown joint distribution $\pi(x, z)$ and the observed current value of $x$, estimate the posterior distribution $\gamma_{pos}(z|x)$.
 
 \item The optimal transport problem (OTP), a particular case of the OTBP with only two marginals, yields a natural horizontal distance among distributions. It also serves as a powerful tool for density estimation and sampling.
 
\end{enumerate}

This article proposes an efficient methodology to solve the OTBP, providing the capability to both simulate and estimate $\rho(x|z)$. It also clarifies the relation between the statistical formulation of the Monge OTBP problem as posed above and the geometrically based Wasserstein barycenter problem \cite{AguehBarycenters}: while the latter addresses the barycenter's distribution $\mu$, the former is centered on the underlying random variable $y$ and its relation to $x$ and $z$. Corollaries include new methodologies for the solution of Monge's OTP and for model-free Bayesian inference.

The methodology's central component is a new adversarial formulation of the pushforward condition, posed in terms of the independence between the random variables $z$ and $y$ and enforced through test functions. When these test functions are restricted to inner-product spaces, the optimal adversarial strategy can be expressed in terms of the first principal components of a matrix, reducing the problem to a pure minimization, which can be solved efficiently through a flow-based gradient-descent procedure. The corresponding optimal map $y = T(x, z)$ can be inverted in closed-form, which facilitates conditional density estimation and simulation. The closed form inversion formula extracts from the data natural factors $\{f^k(z)\}$ that encode the dependence of $x$ on $z$.

\subsection{Relation to prior work}

The Kantorovich --or Wasserstein-- OTBP was introduced in \cite{AguehBarycenters}, defining the barycenter $\mu_*$ of a set of distributions $\{\mu_i\}$ as the minimizer of a weighted sum of the squared Wasserstein distances between the $\{\mu_i\}$ and $\mu_*$. 
As discussed in section \ref{sec:MongeFormulation}, the distribution underlying the solution $y$ to Monge's OTBP agrees with the solution to the Wasserstein barycenter problem extended to general factors $z$, when the latter is supported on $z$-dependent maps. The extension of the OTBP to a continuous covariate $z$ was studied in \cite{pass2013optimal} in the context of its connection to the multimarginal optimal transport in the limit of infinitely many marginals. 

There is a rich literature on the numerical solution of the OTBP, typically in their Kantorovich formulation. 
While discrete methods regard the marginals as convex linear combinations of Dirac delta functions, continuous methods assume that smooth probability density functions underly the data.
The first category includes methods that leverage the Sinkhorn algorithm \cite{essid2019traversing, Sin64, Cuturi} and linear programming-based methods \cite{auricchio2019computing,oberman2015efficient}. Among the first algorithms to treat the problem in a continuous setting are \cite{kuang2019sample, li2020continuous}, both based on the dual of Kantorovich formulation. 
Most algorithms adopting the continuous apporach parameterize the maps pushing forward each $\rho(x|z)$ to the barycenter via a deep neural network with problem-dependent architecture. An example of this approach in \cite{taghvaei2022optimal} uses Convex Neural Networks to parametrize a potential related to the optimal map. A different approach is the flow-based methodology adopted in  \cite{tabak2022distributional} and inspired by \cite{TV, TaTu, TT1}. Flow-based numerical solvers do not require neural networks or any a priori parametrization of the map. They lead naturally to the adoption of gradient descent methods, more straightforward than saddle point optimizers, whose convergence is harder to characterize \cite{mazumdar2019finding, essid2023implicit}.
Another connected work concerns not the OTBP but the related vector quantile regression \cite{Quantiles}. This also transforms a random variable through a factor dependent map to another that is independent of the factor. Unlike the barycenter, however, the target distribution --the quantile-- is fixed by the user, and the transformation proposed is linear and computed through linear programming. 

Our approach focuses on the statistical analysis of data with factors $z$ that typically include continuous components, therefore requiring the solution of OTBP problems with infinitely many marginals. 
A major point distinguishing the work presented here from the existing literature is the interpretation of the pushforward condition that drives the samples underlying the marginals towards the barycenter. As in \cite{tabak2022distributional}, we characterize the pushforward condition in terms of the statistical independence between two random variables: the cofactors $z$ and $y = T(x, z)$. This statistical characterization of the barycenter is critical for a number of applications, such as removal of variability and factor discovery \cite{TT2, yang2020conditional}, and treatment effect estimation \cite{tabak2021data, tabak2024conditional} (See \cite{pmlr-vR5-gretton05a} for a similar characterization of independence through reproducing kernel Hilbert spaces.)
Directly related to the statistical interpretation of the push-forward condition is the ability to solve numerically the barycenter problem for continuous factors and under costs different from the canonical Euclidean distance. While there is some literature dealing with more general costs (see for instance \cite{staib2017parallel,peyre2016gromov}), to the best of our knowledge, the only alternative work on the solution of this problem in the continuous setting --i.e. with infinitely many marginals-- is our own previous work in \cite{tabak2020conditional, yang2020conditional, tabak2022distributional}. The formulations in those articles differ substantially from the current one: the first solved a minimax problem for a potential $\psi(y, z)$, extending the attributable component methodology of \cite{TT3} beyond nonlinear regression, the second developed BaryNet, a network-based algorithm, and the third formulated the push-forward condition in terms of a test function of the form $F(y, z) = \rho(y|z)$, whose estimation through kernels has a computational cost that grows quadratically with the number of samples. By contrast, the current proposal, which is based on flows in phase space and formulates the pushforward condition in terms of the first principal components of a matrix, has a complexity that scales linearly with the number of observations.

Next, we would like to mention existing work involving the specific use of optimal transport in the context of generative modeling. While to the best of our knowledge the first work to use optimal transport for conditional simulation is \cite{tabak2020conditional}, since then a few important contributions have been made \cite{baptista2024conditional,hosseini2025conditional}. Both of these works constrain the map to be triangular and compute the pushforward from an initial density $\rho(x)$ to a conditional density $\mu(x|y)$ given samples from the joint $\pi(x, y)$ distribution. The work in \cite{baptista2024conditional} leverages the Sinkhorn iteration in conjunction with an averaging procedure. The goal is to  output a deterministic pushforward map that, otherwise, would be only a coupling \cite{Cuturi}. This step is crucial for the generation of new sample distributed according to the conditional densities that one wished to simulate. While the use of the Sinkhorn iteration ensures computational efficiency, the averaging procedure introduces an undesirable bias on the target distribution. The work in \cite{hosseini2025conditional}, focuses largely on the well-posedness analysis of the constrained (triangular) problem in the infinite dimensional setting in which the measures to be transported are defined on functional spaces rather than on a final dimensional space. When faced with a practical application, however, these spaces need to be discretized and an efficient algorithm for the solution of the (finite) dimensional optimization problem is still required. For this, the work in \cite{hosseini2025conditional} introduces either a penalty method to minimize the statistical divergence between the starting and the target measures or an adversarial saddle point solver.

The penalty method is very similar to the work on conditional optimal transport first presented in \cite{tabak2021data}, though providing a solid ground to conditional optimal transport in infinite dimensional spaces.
  The solution of the adversarial problem, instead, presents well known convergence issues \cite{essid2023implicit}, which makes the tuning of the solver not always a straightforward task. In this sense, the use of triangular maps has been explored in \cite{alfonso2023generative} in the context of a purely gradient descent flow \cite{TT1} that does not require alternating optimization steps characterizing an adversarial approach.

Besides considerations on the nature of the numerical solvers, the work presented here differs from the conditional optimal transport approach in at least three main conceptual aspects. The first and most important is related to the nature of the barycenter. This distribution is not prescribed a priori like in the case of conditional optimal transport and it constitutes an important byproduct of the algorithm that is functional to a whole series of tasks such as factor discovery and removal of variability (\cite{yang2020conditional}). The second aspect is related to the fact that mapping the marginals $\rho(x|z)$ to the barycenter $\mu(y)$ should, in principle, be cheaper than mapping them to a reference distribution, the reason being that $\mu(y)$ is the closest distribution to the $\rho(x|z)$, while the reference distribution maybe far away from them in the Wasserstein sense. The third aspect is related to the selection of the functional space in which the map is approximated. Our methodology embeds such selection in the numerical solver, a feature that is unique to this approach. Even if one has to choose a priori the broad family of functions that can be used to approximate the optimal map, the algorithm tailors automatically the subset of them that are more relevant for the approximation.

A concluding remark in this section regards the use of transport measure in the area of causal inference \cite{Causal1,Causal2}. This is an important area of research and we leave for future work the study of the barycenter problem within this context.

\subsection{Plan of the article} 

This article is structured as follows. Section \ref{sec:MongeFormulation} discusses the formulation of the Monge OTBP, relates it to an extension of the Wasserstein barycenter problem, and justifies its use to identify hidden sources of variability. Section \ref{sec:Indep} introduces an adversarial formulation of independence, relaxed to finite dimensional inner-product functional spaces. This gives rise to a compact formulation of the problem in terms of the singular values of a matrix, discussed in Section \ref{sec:Flow}, and an efficient procedure for its minimization through gradient descent. Section \ref{sec:mapinv} derives a closed-form expression for the inverse $x = X(y, z)$ of the map $y = T(x, z)$, mediated by extracted factors. It also extends the procedure to pairwise cost functions different from the squared norm. Section \ref{sec:Impl} discusses various implementation aspects: the choice of functional spaces, the determination of the penalization parameter and of the learning rate, the termination criteria and the solution of various barycenter problems in a row, each contributing to further explain the variability of the outcome $x$. It also discusses the algorithm's complexity. Section \ref{sec:OTMonge} solves the regular [Monge] OTP through a suitable reduction of the more general OTBP. This is used in Section \ref{sec:CDE} to perform regular and conditional density estimation. Section \ref{sec:NumEx} illustrates the methodology through numerical examples. We first use synthetic data sets to demonstrate various aspects of the OTBP, such as the simulation and estimation of conditional densities, model-free Bayesian inference and the uncovering of hidden explanatory factors. Then we apply the procedure to real data sets related to weather and climate.  Finally, Section \ref{sec:Summary} summarizes the procedure and suggests avenues for further development. 

\section{A Monge formulation of the optimal transport barycenter problem}\label{sec:MongeFormulation}

Given a joint distribution
$ \pi(x, z) = \rho(x |z)\ \gamma(z) $
between two sets of variables: the \emph{outcome} $x \in \mathcal{X}$ and the \emph{covariates} $z \in \mathcal{Z}$,
we seek a map
$ y = T(x, z) \in \mathcal{X} $
that removes from $x$ the variability that $z$ can explain, i.e. such that the random variables $y$ and $z$ are independent. 
We require that the space $\mathcal{X}$ have the structure of a smooth manifold, while the space $\mathcal{Z}$ can include both continuous and discrete components. We will further assume that $\rho(x|z)$ is absolutely continuous for all $z$, vanishing on small subsets of $\mathcal{X}$.
In order to remove from $x$ only the variability that $z$ can account for, we select the map $T$ that minimizes the expected value of the total transportation cost $ C = \EE_{\pi}[c(x, y)]$,
where $c(x, y)$, an externally provided pairwise cost function, measures the deformation of the data incurred by moving $x$ to $y$. The canonical choice for $c$ is the squared  distance
\begin{equation}
 c(x, y) = \frac{1}{2} \|y - x\|^2.
 \label{canonical_cost}
\end{equation}
More general costs $C$, not necessarily based on pairwise cost functions, give rise to the Distributional Barycenter Problem \cite{tabak2022distributional}. \esteban{In this article, we restrict attention for concreteness to the pairwise canonical cost in (\ref{canonical_cost}), with the slight modifications required to handle more general costs discussed in subsection \ref{General_costs}.}

At first sight, the resulting problem (\ref{Barycenter}) looks quite different from the Wasserstein barycenter problem, introduced in \cite{AguehBarycenters}, which reads
\begin{equation}
\mu_* = \arg \inf_{\mu} \sum_{i=1}^p \lambda_i {W_2}^2\left(\mu_i, \mu\right), \quad
 {W_2}^2(\rho, \mu) = \inf_{\xi(x,y) \in \Pi(\rho, \mu)} \EE_{\xi} \left[\|y-x\|^2 \right] ,
\label{ClassicBarycenter}
\end{equation}
where $\Pi(\rho, \mu)$ is the set of joint distributions having $\rho$ and $\mu$ as marginals.
The differences between (\ref{ClassicBarycenter}) and (\ref{Barycenter}) stem from their conceptual origin: while (\ref{ClassicBarycenter}) extends the geometrical notion of barycenter to sets of distributions equipped with the Wasserstein distance, (\ref{Barycenter}) uses a map $T$ to remove from $x$ any variability that $z$ can explain. 
Yet the two problems relate as follows. 
The $i \in [1, \ldots, p]$ in (\ref{ClassicBarycenter}) correspond to the $z \in \mathcal{Z}$ in (\ref{Barycenter}), their weights $\{\lambda_i\}$ to the distribution $\gamma(z)$, and the $\{\mu_i\}$ to the $\rho(x|z)$. Thus the barycenter problem in (\ref{ClassicBarycenter}) is restricted to discrete covariates $z$, which play the role of indexes for the distributions $\{\mu_i\}$. By contrast, the $z$ in (\ref{Barycenter}) is a random variable of general type linked to $x$ through their joint distribution $\pi$. The well posedness of the barycenter problem with infinitely many marginals has been studied in \cite{pass2013optimal} in connection to the multi-marginal OTP.
Extending Kantorovich's relaxation of the OTP \cite{santambrogio2015optimal}, (\ref{ClassicBarycenter}) considers general couplings $\xi$ between the $\mu_i$ and $\mu$, while (\ref{Barycenter}) extends Monge's original formulation \cite{monge1781memoire} to the barycenter problem, restricting attention to maps $T$ that push forward the $\rho(x|z)$ to $\mu$. These maps are central to the applications that motivate (\ref{Barycenter}), as they are used both for conditional density simulation and estimation. Importantly, they turn $y = T(x, z)$ into a random variable that derives from $x$ and $z$, which leads to another critical distinction: while the argument of the minimization in (\ref{ClassicBarycenter}) is the barycenter $\mu_*$ of the $\mu_i$, the formulation in (\ref{Barycenter}) does not involve the barycenter at all. The fact that $y$ in this formulation is a random variable gives meaning to the alternative requirement of independence between $y$ and $z$. 

That the two problems are much closer than they appear at first sight follows from the fact that, as proved in \cite{Brenier91}, the solution to Kantorovich's formulation of the OTP  for smooth distributions and under quite general assumptions, also solves Monge's, a results that has been extended to the OTBP in \cite{AguehBarycenters,pass2013optimal}. 
We put together the connection between the two problems in the framework of this article through the discussion that follows. First,  the following lemma shows that the two characterizations of the barycenter problem: one geometrical, in terms of the distribution with minimal Wasserstein distance from the marginals, and the other statistical, in terms of the independence between two random variables, are equivalent.

\medskip

\begin{lemma}
Given a joint distribution $\pi(x, z)$, $x \in \mathcal{X}$, $z \in \mathcal{Z}$, define the marginal $ \gamma(z) = \pi(\mathcal{X}, z)$ and the conditional distribution $\rho(x|z) = \frac{\pi(x, z)}{\gamma(z)}$,
and consider the following two problems:
\begin{enumerate}

\item Extended Wasserstein barycenter:
$$ \mu^*, \xi_z^* = \arg\min_{\mu, \xi_z} C_K = \EE_{\gamma} \Big[\EE_{\xi_z}[c(x, y)]\Big],
\quad  \xi_z \in \Pi\left(\rho(x|z), \mu(y) \right) $$
(We call this problem ``extended'' because the covariates $z$ are not necessarily discrete),

\item Monge barycenter:
$$  T^* = \arg\min_{T} C_M = \EE_{\pi}[c(x, y)], \quad y = T(x, z),  \quad y \indep z .$$

\end{enumerate}
If the minimizing couplings $\xi_z^*$ for the first problem are supported on maps,  i.e.,
$ \xi_z^*(x, y) = \pi(x, \mathcal{Z})\ \delta\left(y - Q_z(x)\right)$,
then
$ T^*(x, z) = Q_z(x)$, and consequently $C_M = C_K$ and $\forall z \ \mu^* = T^*(:, z) \# \rho(:|z)$.

\end{lemma}

\begin{proof} 
Since $\xi_z^*$ solves problem 1, $Q_z \# \rho(:|z) = \mu^*$, so the joint distribution satisfies
$ \Theta(y, z) = \mu^*(y) \gamma(z)$,
which implies that $y \indep z$. If there existed another $T \ne Q_z$ pushing forward all $\rho(x|z)$ to a single distribution $\mu(y)$ at a cost $C_M < C_K$, then $ \xi_z(x, y) = \pi(x, \mathcal{Z})\ \delta\left(y - T(x, z)\right)$ would solve problem 1 with a smaller cost than the optimal $\xi_z^*$, a contradiction.

\end{proof}

If the joint distribution $\pi(x, z)$ is absolutely continuous in $x$ for a set of $z \in Z$ of nonzero measure, then the barycenter $\mu(y)$ of the corresponding $\rho(x|z)$ under the canonical cost in (\ref{canonical_cost})) is also absolutely continuous and the optimal couplings $\xi_z^*$ are indeed supported on maps. 
This result was proved for discrete $z$'s in \cite{AguehBarycenters},  Theorem 5.1, and extended to continuous $z$  in \cite{pass2013optimal}, see Corollary 3.3.3 and Theorem 4.2.5.
We might conclude that the Monge and [extended] Wasserstein barycenter problems are equivalent when applied to smooth distributions: after all, their unique solutions map to each other. Yet this equivalence applies to the problems' solutions, not to their formulations. All the applications described in this article, as well as the methodology proposed for solving the problem numerically, are strongly based on the map $y = T(x, z)$ and its inverse $X(y, z)$, both parameterized by $z$ and defining the random variables $(x, y)$ in terms of each other.

The relation between a random variable $x \in \mathcal{X}$ and known covariates $z \in \mathcal{Z}$ can be specified alternatively through the conditional distribution $\rho(x|z)$ and through a functional relation
$ x = \phi(w, z)$, $w \in \mathcal{W}$, $w \sim \nu(: | z)$,
where the random variable $w$ represents all additional causes of variability in $x$, which we either do not currently consider, cannot measure or are simply not aware of. The function $\phi$ and the distribution $\nu$ underlying $w$ determine $\rho(x|z)$ uniquely, but more than one pair $(\nu, \phi)$ can give rise to the same $\rho$. One special pair is provided by the solution to the OTBP:

\begin{lemma}\label{th:hidden}
Given any joint distribution $\pi(x, z)$ that vanishes on small sets in $\mathcal{X}$ for all $z$, the random variable $x$ can be written as
$ x = X(y, z)$,
where $y = T(x, z) \sim \mu$ is the solution to the barycenter problem for $\rho(x|z)$, so $X(:, z) = T^{-1}(:, z)$.
\end{lemma}

\begin{proof} 
Since neither $\rho(x | z)$ (for any $z$) nor $\mu(y)$ assign finite measure to small sets, the fact that the least costly map $T(:, z)$ pushing forward $\rho(: | z)$ to $\mu$ is invertible is a central result in optimal transport theory \cite{santambrogio2015optimal}. 

\end{proof}

This lemma provides the ground for various applications. 
\begin{enumerate}
\item In order to generate samples $\{x_i^*\}$ from $\rho(x | z_*)$ for any target value $z_*$, it is enough to generate samples $\{y_i\}$ from the barycenter $\mu(y)$ and write $x_i^* = X(y_i, z_*)$. A number $n$ of such samples is already available through the barycentric map acting on the available data pairs, $y_i = T(x_i, z_i)$. More can be obtained if needed, performing optimal transport between $\mu$ and a given distribution $\nu$ and pushing forward samples of the latter to $\mu$ through the inverse of the corresponding transportation map. This is detailed in sections \ref{sec:OTMonge} and \ref{sec:CDE}.

\item Since our algorithm provides both $X(y, z)$ and its $y$-derivatives (see section \ref{sec:mapinv}), one can estimate $\rho(x | z_*)$ by the change-of-variable formula applied to an estimate for $\mu(y)$ (which can itself be obtained by optimal transporting $\mu$ to a known $\nu$ and applying again the change-of-variable formula.)

\end{enumerate}

A third application addresses the following question: assuming that there exists a ``true'' additional source $w$ of variability in $x$, such that  that $x = \phi(w, z)$, how much does our $y = T(x, z)$ teach us about the true $w$? (Here the notion of a true source is field dependent; for our purposes, we just assume that such true $w$ exists.)
An identifiability issue arises: without additional information, the conditional distribution $\rho(x|z)$ does not suffice to determine $w$. For instance, if for some value of $z$, $\phi\left(w_1, z\right) = \phi\left(w_2, z\right)$, then there is no way that using $x$ and $z$ alone we could distinguish between $w_1$ and $w_2$. More generally, any two $z$-dependent random variables $W_1^z$ and $W_2^z$ such that the distributions of both $\phi_1(W_1^z, z)$ and $\phi_2(W_2^z, z)$ agree with $\rho(x|z)$ explain the data equally well. In particular, $y$ provides one such explanatory variable, with the additional property that it is necessarily a function of $w$ and $z$:
$ y = T(x, z) = T(\phi(w, z), z) = Y(w, z)$.
Moreover, $Y(:, z)$ is invertible for all values of $z$ for which $w$ is identifiable, i.e. such that $\phi\left(w_1, z\right) = \phi\left(w_2, z\right) \Rightarrow w_1 = w_2$, since
$$ Y\left(w_1, z\right) = Y\left(w_2, z\right)\ \Rightarrow \ T\left(\phi(w_1, z), z \right) = T\left(\phi(w_2, z), z \right) \ \Rightarrow \ \phi(w_1, z) = \phi(w_2, z) \ \Rightarrow \ w_1 = w_2. $$
Then, in order to uncover the ``true'' hidden explanatory $w$, one can use the fact that its identifiable component must have a --possibly $z$-dependent-- one-to-one relation to $y$, together with any other information available on $w$, such as other variables that it may depend upon or correlate with. In the absence of any such additional information, $y$ is the most natural explanatory variable among all $w$, since by construction it is independent of the known factors $z$, it is the closest to $x$ itself, and it is the most ``economical'', since it is identifiable for all $z$. We illustrate these concepts through examples in Section \ref{sec:NumEx}.

\section{Adversarial characterization of independence}\label{sec:Indep}

Posing the Monge OTBP (\ref{Barycenter})
in data-driven scenarios, where the joint distribution $\pi(x, z)$ is only known through $n$ sample pairs $\{x_i, z_i\}$, requires a sample-friendly formulation of the independence condition between the random variables $y$ and $z$. We will use a weak characterization based on test functions \cite{jacod2012probability,tabak2022distributional}:
two variables $y \in \mathcal{Y}$ and $z \in \mathcal{Z}$ with joint distribution $\pi(y, z)$ are independent if and only if any two bounded measurable functions $g(y)$ and $f(z)$ satisfy
$$ \EE_{\pi}[g(y) f(z)] = \EE_{\rho}[g(y)]\ \EE_{\gamma}[f(z)], \quad
 \rho(y) \defeq \pi(y, \mathcal{Z}), \ \gamma(z) \defeq \pi(\mathcal{Y}, z). $$
We can restrict the functions $f$ and $g$ to smaller spaces $\mathcal{F}$ and $\mathcal{G}$, provided that \esteban{these} contain suitable approximations to the delta function. \esteban{This follows from the fact} that,  if $ \pi\left(y, z\right) \ne \rho\left(y\right) \gamma\left(z\right)$ in a neighborhood $U_*$ of the pair $(y_*, z_*)$, 
then $ g(y) = \delta_{y_*}\left(y\right)$, $f(z) = \delta_{z_*}\left(z\right)$
satisfy
$ \EE_{\pi}[g(y) f(z)] \ne  \EE_{\rho}[g(y)]\ \EE_{\gamma}[f(z)]$,  
where $\delta_{y_*}$ and $\delta_{z_*}$ are non-negative functions whose product is positive at $(y_*, z_*)$ and vanishes outside of $U_*$. \esteban{Then  the spaces $\mathcal F$ and $\mathcal G$ need not include functions that vary at scales smaller than those that the available data allow us to resolve. We adopt for  $\mathcal{F}$ and $\mathcal{G}$ finite dimensional inner product linear spaces, i.e. spaces spanned by a finite set of functions, provided with the canonical inner product over the available sample points.}

Subtracting the mean of $f(z)$ yields the following equivalence statement:
two variables $y$ and $z$ are independent if and only if, for all functions $g(y) \in \mathcal{G}$ and $f(z) \in \mathcal{F}$ with $\EE_{\gamma}[f(z)] = 0$,
$ \EE_{\pi}[g(y) f(z)] = 0$.
This equivalence gives rise to the following adversarial formulation of the barycenter problem (\ref{Barycenter}):

\begin{equation}
 \min_{y = T(x, z)} \max_{g, f, \lambda} L = \EE_{\pi}[c(x, y)] + \lambda\ \EE_{\pi}[g(y) f(z) ], \quad \EE_{\gamma}[f] = 0,\ \|f\| = \|g\| = 1,
 \label{Barycenter_adversarial_lambda}
\end{equation}
where we have decoupled the amplitude of $f$ and $g$ from their shape, absorbing their amplitude in the factor $\lambda$. Moreover, we can replace the maximization over $\lambda$ by the external provision of a a penalization parameter  $\lambda \gg 1$ for non-compliance of the independence condition,  a relaxation that converges to (\ref{Barycenter_adversarial_lambda}) as $\lambda \rightarrow \infty$.

\esteban{
The norms of $f$ and $g$ defined through the canonical inner products
\begin{equation}
\left(f_1, f_2\right) = \EE\left[f_1, f_2\right], \quad \left(g_1, g_2\right) = \EE\left[g_1, g_2\right], 
 \label{inner_product_1}
\end{equation}
 represent the standard deviation of $f$ and $g$, since not only $\EE_{\gamma}[f] = 0$ by construction, but also $\EE_{\rho}[g] = 0$ holds at the optimal solution: a constant added to $g$ does not affect the value of $L$, and the norm of $g - a$ is smallest when $a = \bar{g}$. It follows that we could read the problem as the minimization of the transportation cost subject to the condition that the correlation between any two functions $f(z)$ and $g(y)$ vanishes \cite{jacod2012probability}. 
}

A data-driven formulation of (\ref{Barycenter_adversarial_lambda}) replaces expected values by empirical means,
\esteban{
\begin{equation}
 \min_{y_i = T\left(x_i, z_i\right)} \max_{g, f} L = \sum_{i=1}^n \left[\frac{1}{n}  c\left(x_i, y_i\right) + \lambda\ g\left(y_i\right) f\left(z_i\right)\right], \quad \sum_{i=1}^n f\left(z_i\right) = 0,\ \|f\| = \|g\| = 1.
 \label{Barycenter_sample_based}
\end{equation}
}
Since we cannot enforce infinitely many constraints on the finite set $\{y_i\}$ without trivializing the solution,
we supplement (\ref{Barycenter_sample_based}) with the specification of two finite dimensional inner-product spaces of functions $\mathcal{F}$ and $\mathcal{G}$ over which to perform the maximization, writing
$$  f(z) = F(z) a, \quad  g(y) = G(y) b, \quad a \in R^{m_z}, \quad  b \in R^{m_y}, $$
where the $m_z$ columns of $F$ and the $m_y$ columns of $G$ are functions respectively of $z$ and $y$ (both have $n$ rows when evaluated at the sample points) and the functions acting as columns of $F$ have zero mean.
Choices for the functions defining $\mathcal{F}$ and $\mathcal{G}$ will be discussed in Section \ref{sec:FandG}.
Independently of their choice, some further processing is required, which we describe here in terms of $\mathcal{G}$, since the same process applies to $\mathcal{F}$.  

\esteban{We replace $G$ --an operator from $R^{m_y}$ to $\mathcal{G}$-- by an orthogonal operator $Q_y$ with the same effective range, so that the requirement that $\|g(y)\| =1$ translates into the condition that $\|b\| = 1$. 
To this end, we divide the procedure into stages. Calling $\{y_j^0\}$ the values of $\{y_j\}$ at the beginning of a stage,
we perform the reduced singular-value decomposition
$ G_i^j \defeq G^j\left(y^0_i\right) \approx \sum_{k=1}^{n_y} \sigma_k\ u^k_i \ v^k_j$,
with $n_y \le m_y$ chosen that the sum $ \sum_{k=1}^{n_y} \left(\sigma_{k}\right)^2 $
is larger that a fraction of the squared norm of $ \left\|G\right\|^2 = \sum_{i,j} \left(G_i^j \right)^2 = \sum_{k = 1}^{m_y} \left(\sigma_k\right)^2$, and we adopt 
$ Q_{y}(y) = G(y) B^y, \quad B^y_{jk} = \frac{1}{\sigma_k} v^k_j$
(Notice that, in particular, $Q_y^k\left(y^0_i\right) = u^k_i$.)
As the $\{y_j\}$ evolve through the procedure, we need to update $Q_y$ (and correspondingly $B$) so that $Q_y$ remains orthogonal, which for a rectangular matrix reads $Q_y' Q_y = I$. Rather than applying Gram-Schmidt at each step, which would be costly, we make use of the following lemma:
\begin{lemma}\label{th:QtoQ}
If a matrix $Q \in \R^{n\times m}$ satisfies $Q' Q = I_m + O(\epsilon)$, then $\tilde{Q} = Q \beta$, with $\beta = \frac{3 I - Q'Q}{2}$, satisfies $\tilde{Q}' \tilde{Q} = I_m + O\left(\epsilon^2\right)$.
\end{lemma}
\begin{proof} 
Write $Q' Q = I_m + \epsilon C$, with $C = O(1)$. Then $\beta = \frac{3 I - Q'Q}{2} = I_m - \epsilon \frac{C}{2}$, and 
$$ \tilde{Q}' \tilde{Q} = \beta' Q' Q \beta = \left(I_m - \epsilon \frac{C}{2} \right) \left(I_m + \epsilon C \right) \left(I_m - \epsilon \frac{C}{2} \right) = I_m - \frac{3C^2}{4} \epsilon^2 + \frac{C^3}{4} \epsilon^3 . $$
\end{proof}
\noindent
(A complementary lemma with no small parameter $\epsilon$ proves that, under the same notation, if $Q \ne 0$ and the largest eigenvalue $\lambda_{max}$ of $Q' Q - I$ is smaller than $4$, then the largest eigenvalue of $\tilde{Q}' \tilde{Q} - I$ is smaller than $\lambda_{\max}$ and iterating the procedure (i.e. applying it again to $\tilde{Q}$ and so on) converges to an orthogonal matrix.)\\
In view of this lemma, to preserve the orthogonality of $Q$ after each update of $y$, we update $B$ to $B \beta$, where $\beta = \frac{3 I - Q'Q}{2}$, with $Q = G(y) B$.\\
The reason for dividing the procedure into stages is that the range of $G$ depends on the distribution of $y$, which evolves through the procedure. Then the orthogonal $Q_y$ that best captures this range needs to be updated as $y$ changes. Moreover, as discussed in subsections \ref{ChoicesFG} and \ref{sec:FandG}, it is convenient to adapt the functional space $\mathcal G$ to the current $y$ as well. It is much more economical to perform these updates in discrete stages than at every step of gradient descent.
}

\esteban{
Unlike $y$, the covariates $z$ do not evolve during the procedure, so we apply the reduced singular value decomposition to $F(z)$ only once, producing an orthogonal matrix $Q_z$ of rank $n_z$ and corresponding matrix $B^z$.
}

\section{A flow-based methodology}\label{sec:Flow}

Replacing in (\ref{Barycenter_sample_based}) $f(z_i)$ by $Q_z(z_i)  a$ and $g(y_i)$ by $Q_y(y_i) b$ yields
\esteban{
$$
 \min_{\{y_i\}} \left\{\max_{a, b} \frac{1}{n} \sum_{i=1}^n c\left(x_i, y_i\right) + \lambda \sum_{h=1}^{n_z} \sum_{l=1}^{n_y} \left(\sum_{i=1}^n Q_z^h(z_i) Q_y^l(y_i)\right) a_h b_l , \quad \|a\| = \|b\| = 1 \right\}.
 $$
 }
The maximization over $a$ and $b$ can be carried out explicitly: they must align with the left and right first principal components of the $n_z\times n_y$ matrix
$ A^{hl} \defeq \sum_i Q^h_z\left(z_i\right) Q^l_y\left(y_i\right)$,
and the penalty term is given by the first singular value $\sigma_1(y)$ of $A$. 
It follows that we can write the problem as a minimization over $y$ alone: 
\esteban{
\begin{equation}
 \min_{\{y_i\}} L = \frac{1}{n} \sum_{i=1}^n c\left(x_i, y_i\right) + \lambda \left\|A(y)\right\|,  
 \quad \left\|A\right\| \defeq \sigma_1 = \max_{\|a\| = \|b\|=1} a' A b.
 \label{Barycenter_y}
\end{equation}
}
We can interpret the corresponding functions $f(z) = Q_z(z) a$, $g(y) = Q_y(y) b$ as the features whose correlation most strongly displays the current dependence between $z$ and $y$.

This suggests a flow-based procedure, whereby $y$, initially set equal to $x$, follows gradient descent of (\ref{Barycenter_y}),
\begin{equation}
  y_i^{n+1} = y_i^n - \eta^n \left[\frac{1}{n} \nabla_y \left. c\left(x_i, y\right)\right|_{y_i^n}  + \lambda\   a' \left. \nabla_y A \right|_{y_i^n} b  \right], \quad \left. \nabla_y A^{hl} \right|_{y_i^n} = Q^h_z\left(z_i\right) \sum_j \left. \nabla G^{j}(y)\right|_{y_i} B^y_{jl} ,
\label{GD}
\end{equation}
 for which all $\{y_i\}$ decouple,
and $a$ and $b$ are updated in an alternate step.

It might appear that we are taking an uncontrolled approximation to the $y$-gradient of $L$ in (\ref{Barycenter_y}) by differentiating only $A$ in (\ref{GD}) at fixed $a$ and $b$. The principal components of $A$ do of course depend on $A$, so they too change when $A$ varies. Yet this way of computing derivatives is exact:

\medskip

\begin{lemma}\label{th:dsigma}
The derivative of the $j$-th principal value $\sigma_j$ of a matrix $A$ with respect to any parameter $s$ on which $A$ may depend, is given by
$$ \frac{\partial}{\partial s} \sigma_j = a'  \left(\frac{\partial}{\partial s} A \right) b, $$
where $a$ and $b$ are the left and right $j_{\text{th}}$ principal components of $A$.

\end{lemma}

\begin{proof} 
By definition, $ \sigma_j = a' A b$, so
$$ \frac{\partial}{\partial s} \sigma_j = a'  \left(\frac{\partial}{\partial s} A \right) b
+ \left(\frac{\partial}{\partial s} a \right)' A b + a'  A \left(\frac{\partial}{\partial s} b \right). $$
But the principal components satisfy $ A b = \sigma_j a$, $A' a = \sigma_j  b$ and $\|a\| = \|b\| = 1$,  so
$$ \left(\frac{\partial}{\partial s} a \right)' A b = \sigma_j  \left(\frac{\partial}{\partial s} a \right)' a
= \sigma_j \frac{\partial}{\partial s} \frac{\|a\|^2}{2} = 0 \quad
\hbox{and} \quad
a'  A \left(\frac{\partial}{\partial s} b \right) = \sigma_j b' \left(\frac{\partial}{\partial s} b \right) 
= \sigma_j \frac{\partial}{\partial s} \frac{\|b\|^2}{2} = 0 .$$

\end{proof}

The penalty term $\sigma_1(y)$  is not smooth at its arg-min $y = y^*$: for $y$ and $z$ to be independent, $A$ must vanish, and the first singular value $\sigma_1(y)$ of a matrix that depends smoothly on $y$ typically has corners where $A(y)$ vanishes (The simplest example is the $1\times 1$ matrix $A = y \in R$, whose only singular value $\sigma = |y|$ has a corner at $y = 0$.)
To address this, we square the penalty term:
$
 \min_{\{y_i\}} L = \sum_{i=1}^n c\left(x_i, y_i\right) + \lambda \  \sigma_1^2(y)$. 
There still remains an issues to address to make the methodology fully functional. Because every step of the algorithm brings down the largest singular value of $A(y)$, the first few of those singular values will tend to coalesce at convergence at a common value $\sigma_* \ll 1$. Then the derivatives of the penalty term with respect to the $\{y_i\}$ are not well defined, as they depend on the arbitrary choice of one pair among the various singular components $(a, b)$ associated to the singular value $\sigma_*$. In terms of test functions, more that one pair of functions $(f, g)$ have reached the threshold correlation $\sigma_*$.
To address this, we modify the algorithm so that it tracks the first $J$ pairs $(a_j, b_j)$ of principal component of $A$, where $J = \min(\hbox{rank}(A), J_{max})$, with $J_{max}$ fixed by the user.
Then we descend over $y$
\esteban{
\begin{equation}\label{L_multiK_2}
 \min_{\{y_i\}} L = \frac{1}{n} \sum_{i=1}^n c\left(x_i, y_i\right) + \sum_{j=1}^{J} \lambda_j\ {\sigma_j}^2\left(y\right),  \quad
 \sigma_j\left(y\right) \defeq {a_j}' A(y) b_j .
\end{equation}
}
Notice that this extension carries little computational cost, since the $A(y)$ to differentiate is common to all the $\{\sigma_j\}$. If performed using reproducing Kernel Hilbert spaces, this extension could be thought as interpolating between COCO \cite{pmlr-vR5-gretton05a} and HSIC \cite{10.1007/11564089_7}, where it is shown both that the penalty term vanishes if and only if $y$ and $z$ are independent,  and that its empirical estimate converges to the population value at a rate $\frac{1}{\sqrt{n}}$.

\section{Map inversion}\label{sec:mapinv}

The procedure described so far finds $n$ samples $y_i = T(x_i, z_i)$ of the barycenter $\mu(y)$. In order to simulate $\rho(x | z_*)$ for a target $z_*$, we need to invert $T$ to obtain $n$ samples $\{x^*_i\}$ from $\rho(: | z_*)$ through
$$ x^*_i = X\left(y_i, z_*\right) \defeq T^{-1}\left(y_i, z_*\right). $$
Since we do not know $T(x, z)$ in closed form, it could appear that we can only invert it by learning $X(y, z)$ from its $n$ available samples $\{x_i, y_i, z_i\}$, for instance through kernel regression, nearest neighbor or neural networks.

\esteban{
Yet we can do much better than that and obtain a closed form for $T^{-1}$, exploiting the fact that the penalization parameters $\{\lambda_j\}$ are large but finite. Since at convergence the gradient $\nabla_y L$ is zero, or at least sufficiently small to satisfy a termination criterion,  we have
$ \nabla_{y_i}  c\left(x_i, y_i\right) + 2 n 
 \sum_{j=1}^{J} \lambda_j \ \sigma_j \nabla_y \sigma_j \Big|_{y_i, z_i}= 0$,
which for the canonical quadratic cost in (\ref{canonical_cost}) yields $ x_i = y_i + 2 n  \sum_{j=1}^{J} \lambda_j \ \sigma_j \nabla_y \sigma_j \Big|_{y_i, z_i}$.
In order to invert the map for arbitrary values of $y$ and $z$, we extend the validity of this expression and write
\begin{equation}
 X(y, z) = y +  2 n
 \sum_{j=1}^{J} \lambda_j \ \sigma_j \nabla_y \sigma_j \Big|_{y, z},
\quad \hbox{where} \quad \nabla_y \sigma_j \Big|_{y, z} = f_j(z) \nabla g_j(y),
 \label{X(y)}
\end{equation}
an expression that is smooth in $(y, z)$ and yields $X(y^i, z^i) = x^i$ on all the available samples.}

The inversion formula in (\ref{X(y)}) provides us with a valuable bonus: it shows that the dependence of $x$ on $z$ that our algorithm has uncovered is mediated by the $J$ functions $\{f_j(z)\}$, so we have inadvertently performed \emph{factor extraction}. 
These factors bring in insights about the mechanisms of the dependence of $x$ on $z$, while their gradient inform us of the sensitivity of $x$ with respect to changes in $z$.
When we consider density estimation in Section \ref{sec:CDE}, it will be useful to notice that we have access not only to $X(y, z)$ but also to its derivatives,
\esteban{
\begin{equation}
 \frac{\partial X^p(y, z)}{\partial y^q} = \delta_p^q +  2 n
 \sum_{j=1}^{J} \lambda_j \ \sigma_j \frac{\partial^2}{\partial y^p \partial y^q} \sigma_k  \Big|_{y, z},
\quad \frac{\partial^2}{\partial y^p \partial y^q} \sigma_j  \Big|_{y, z} = f_j(z) \left(\frac{\partial^2 G(y)}{\partial y^p \partial y^q}  \right) B^{y} b_j.
 \label{Xy(y)}
\end{equation}
\subsection{More general costs}
\label{General_costs}
Even though this article restricts attention to the canonical quadratic cost $c(x, y)$ in (\ref{canonical_cost}), the default choice in normed spaces, the problem in hand may call for a different cost function. Examples include costs given by the squared geodesic distance in Riemannian manifolds, by an action when a natural Lagrangian function is available, even by costs that depend on the underlying distribution itself, useful for learning tasks \cite{Zichu}. One may also want to select cost functions that are less sensitive to outliers.\\
The procedure developed uses little of the individual structure of the canonical cost, and thus extends naturally to more general cost functions. Solving the OTBP only requires that an expression for the gradient $\nabla_y c(x, y)$ be available for use in the gradient descent of $L$ in (\ref{L_multiK_2}). In order to extend the formula (\ref{X(y)}) for the inversion of $T(x, z)$, we need this gradient to be also invertible for $x$:
$ w = \nabla_y c(x, y) \Rightarrow x = d(w, y)$ (Notice incidentally that this is precisely the condition under which a solution to Kantorovich relaxation of the OT problem yields a Monge solution.) 
Then, in terms of this function $d$, (\ref{X(y)}) extends to
\begin{equation}
 X(y, z) = d\left(- 2 \lambda\ 
 \sum_{j=1}^{J} \sigma_j \nabla_y \sigma_j \Big|_{y, z}, y\right).
 \label{X(y)_gen}
\end{equation}
For instance, for the cost function
$$ c_p(x, y) = \frac{1}{p} \|y - x \|^p, \quad p > 1,$$
a natural extension of the canonical cost $c_2$, we have
$$ w = \nabla_y c_p(x, y) = \|y - x \|^{p-2} (y - x), \quad 
d(w, y) = y - \frac{w}{\|w\|^{\frac{p-2}{p-1}}} . $$
}

\section{Implementation}\label{sec:Impl}

\esteban{
This section discusses various implementation issues: descent methodology, termination criteria, choice of the penalization parameters $\{\lambda_j\}$ and functional spaces $\mathcal F$ and $\mathcal G$, and the algorithm's computational complexity.}

\esteban{
\subsection{A summary of the procedure}
In order to make the general procedure and its various components clearer to the reader, we summarize them here (Note that some of these components are developed further down the article.)
We minimize the objective function $L$ in (\ref{L_multiK_2}) through gradient descent:
$ y_i^{n+1} = y_i^n - \eta^n \left.\nabla_{y_i} L\right|_{y^n}$, 
back-tracking the learning rate $\eta^n$ from $\theta \eta^{n+1}$, $\theta > 1$, through the Armijo-Goldstein rule \cite{goldstein1965steepest,armijo1966minimization}.
A sketch of the main procedure for the OTBP is provided by the pseudo-code in Algorithm 1. \\
\begin{algorithm}
    \begin{algorithmic}[1]
       \State \textbf{Input:} $n$ pairs $\{x_i, z_i\}$
       \State \textbf{Select:} Functional spaces $\mathcal F$, $\mathcal G$ to use, pairwise cost function $c(x, y)$ \Comment{Default cost $c = \frac{1}{2}\|y - x\|^2}$
       \State \textbf{Set:} $\sigma_* \propto \frac{1}{\sqrt{n}}$ \Comment{Target maximal correlation between $g(y)$ and $f(z)$, $g \in \mathcal G$, $f \in \mathcal F$}
\                     \State {\textbf{Initialize:}} State $y = x$ and learning rate $\eta = \eta_0$.
                     \State Compute the matrix $F_{ij} = F^j\left(z_i\right)$ and $Q_z = F B_z$, $Q_z' Q_z = I$.
       \For{$j_s = 1$ \textbf{to} $n_s$} \Comment{Outer loop over stages.}
       \State Set $y^0 = y$.
       \State Compute the current $G_{ij} = G^j\left(y_i\right)$ and $Q_y = G B_y$. \Comment{Also $F$ if $\mathcal F$ enriched at each stage}
\While{not converged or done with stage $j_s$} \Comment{Inner loop of gradient descent of $L$ over $y$}
\State Compute $L$, $\nabla_y L$, the principal components $\{a_j, b_j\}$ of $A = Q_z' Q_y$ and the $\{\lambda_j\}$.
\State Determine learning rate $\eta$ using backtracking line search satisfying the Armijo condition.
\State Update $y \leftarrow y - \eta \nabla_y L$.
\State Update $B_y \rightarrow B_y \beta$, $Q_y = G(y) B_y$. \Comment{Lemma \ref{th:QtoQ}}
\State Done with stage $= \|y - y^0\|^2 > \frac{1}{10} \|y^0 - \overline{y^0}\|^2$. \Comment{Time to update $y^0$, $\mathcal G$ and $B_y$.}
\State Converged $= \forall j\ \sigma_j < \sigma_*\ \& \ \|\nabla_y L\|^2 < $ threshold.
\EndWhile
\EndFor
\State \textbf{Return:} Final target samples $y = (y_1, \dots, y_n)$, $y_i = T(x_i, z_i)$, $y \indep z$.
    \end{algorithmic}
    \caption{Main OTBP algorithm}
     \label{alg:main_algorithm}
\end{algorithm}
Once OTBP has been run, one can draw $n$ independent samples from $\rho(x|z_*)$ through formula (\ref{X(y)}), where $f_j(z_*) = F(z_*) B_y a_j$. In order to draw more than $n$ samples and/or to estimate $\rho(x|z_*)$, not just simulate it, one applies the procedure of section \ref{sec:CDE}. This further transports the barycenter $\mu$ to a known distribution $\nu$, which is used both to estimate and sample $\mu(y)$, from which $\rho(x|z_*)$ can be estimated and sampled through $T^{-1}(y, z_*)$.
}

\subsection{Termination criterion}

At convergence, $y = y^*$ must satisfy two natural criteria for termination:
\begin{enumerate}

\item Any remaining dependence of $y$ on $z$ must be within an acceptable range: 
\begin{equation}
 \forall j \ \sigma_j\left(y^*\right) < \sigma_* \ll 1. 
 \label{term1}
\end{equation}
In order to assign a value to $\sigma_*$, notice that $\sigma_j$ represents the empirical correlation between $f_j(z)$ and $g_j(y^*)$, which should be uncorrelated for all $f \in \mathcal{F}$, $g \in \mathcal{G}$. 
It follows that a reference value for $\sigma_*$ is the standard deviation of the empirical correlation between two independent variables, which equals $\frac{1}{\sqrt{n}}$, so we have adopted for our experiments
\begin{equation}
  \sigma_* = \frac{0.2}{\sqrt{n}}.
\end{equation}

\item The gradient $\nabla_y L$ of the objective function must be sufficiently small for the inversion formula (\ref{X(y)}) to be valid. Since the error in the determination of $x_i$ from this formula is given by
$$ \left\|x_i - X(y_i, z_i)\right\| = \left\| \nabla_{y_i} L\big|_{y^*}\right\|$$
and a natural reference scale for the square of this error is the variance of $x$, we use as termination criterion
\esteban{
\begin{equation}
\sum_i \left\| \nabla_{y_i} L\big|_{y^*}\right\|^2 < \alpha\ \hbox{var}(x),
 \label{term2}
\end{equation}
}
with $\alpha \ll 1$. We end the run when both (\ref{term1}) and (\ref{term2}) are satisfied.

\end{enumerate}

Another, internal termination criterion starts a new stage once the current $y$ differ significantly from their values $y^0$ at the outset of the current stage, i.e. when
$ \sum_{i=1}^n \left\|y_i - y_i^0\right\|^2 >  \delta \sum_{i=1}^n \left\|y^0_i - \bar{y^0}\right\|^2, \quad 0 < \delta < 1$.
We have adopted in our experiments $\alpha = 0.0025$ and $\delta = 0.1$. 

\esteban{
\subsection{Choice of the penalization parameters}
The penalization parameters $\{\lambda_j\}$ establish a balance at the final $y = y^*$ between the gradients of the transportation cost function and the penalization terms $\{\sigma_j^2\}$. 
We can evolve the $\{\lambda_j\}$ based on the target maximal correlation $\sigma_*$ between the $g$'s and $f$'s. We have from (\ref{L_multiK_2})  that
$$ \nabla_{y_i} L = \frac{y_i - x_i}{n} +  2 \sum_{j=1}^K \lambda_j\ \sigma_j \nabla_{y_i} \sigma_j, $$
where the $\{|\sigma_j|\}$ are sorted in decreasing order.  Then we determine the corresponding $\{\lambda_j\}$ in reverse sequential order (i.e. from $j=J$ to $j=1$) as follows. Assume that we have already determined all $\{\lambda_j\}$ from $j = l+1$ to $J$. Then, compute
$$ u_i = \frac{y_i - x_i}{n} +  2 \sum_{j=l+1}^{J} \lambda_j \ \sigma_j \nabla_{y_i} \sigma^j , \quad v_i = 2 \sigma_l \ \nabla_{y_i} \sigma_l.$$ 
Here $u \in \R^{n\times d_x}$ contains all components of $\nabla_{y} L$ already known, and $v$ contains the component of $\nabla_{y} L$ along the direction of the $l$'th penalty, except for the factor $\lambda_l$ that we seek to determine. We choose $\lambda_l$ so that, if $\sigma_l > \sigma_*$, $\nabla_y L$ with all penalty components from the $l$'th component on projects positively onto $v$, so that $\sigma_l$ decreases in the direction of the gradient of $L$:
$$  \left(\left(\frac{y_i - x_i}{n} +  2 \sum_{j=l}^J \lambda_j\ \sigma_j \nabla_{y_i} \sigma_j\right) , v\right) = (u, v) + \lambda_l \ \|v\|^2 \ge \gamma \|v\|^2,  
\quad \gamma > 0, $$
or
$$ c_{uv} + \lambda_l \ge \gamma, \quad  c_{uv} = \frac{(u, v)}{\|v\|^2} = \frac{\sum_i u_i \cdot v_i}{\sum_i \|v_i\|^2}. $$
We adopt
$$ \gamma = \min \left(\max\left(\frac{\sigma_l}{\sigma_*}  - q, 0\right), \gamma_{max}\right), \quad 0 < q \le 1, \quad \gamma_{max} > 0$$
and set
$$ \lambda_l =\gamma-\min\left(c_{uv},0\right). $$ 
Adopting a threshold  $q$ smaller than $1$ guarantees that all $\sigma_l$ will reach values below $\sigma_*$, enforcing the constraint that $\forall j \ \sigma_j \le \sigma_*$. For our experiments, we have adopted $q=0.9$ and $\gamma_{max} = 20$. \\
With the $\lambda$'s evolving adaptively, the definition of the objective function $L$ changes from step to step, making it hard to guaranty that gradient descent will converge, in the sense of satisfying the second convergence criterion of the prior subsection. To address this, after  (\ref{term1}) is satisfied,  we freeze the $\{\lambda_l\}$ at their current values.
}

\subsection{Choices for the functional spaces $\mathcal{F}$ and $\mathcal{G}$}
\label{ChoicesFG}

\esteban{In order to fully specified the methodology, we need to select the functional spaces $\mathcal{F}$ and $\mathcal{G}$, defined respectively by the linear span of features $\{F^i(z)\}$ and $\{G^j(y)\}$. We discuss in this subsection some general principles guiding their choice, and in an appendix the specific choices that we have made for the experiments in this article. We note though that the choice of features is a broad problem permeating data science; a thorough exploration of optimal choices within the framework of our OTBP would vastly exceed this article's scope.  \\
Since our descent procedure is divided into stages, we can upgrade the functional spaces at the onset of each stage. This is more significant for $\mathcal G$ than for $\mathcal F$, since the $\{z_i\}$ are fixed, while by the end of each stage, the $\{y_i\}$ have evolved so that all current features $G^j(y)$ are independent of $z$ (as captured by $\mathcal F$.)
}

\esteban{
When $X = R^d$, the simplest choice for $\mathcal{G}$ is the space of linear functions,  enforcing independence between any $f(z)$, $f \in \mathcal{F}$ and the conditional mean of $y$. Then, from (\ref{X(y)}),
$ X(y, z) = y + h(z)$, where $ h(z) = 2
\sum_{j=1}^{J} \lambda_j \sigma_j\ f_j(z)\ v_j$  ($v_j \defeq \nabla_y g_j(y)$ is a constant, since $g_j(y)$ is linear.)
It follows that $ h(z) = \bar{x}(z) - \bar{y}$,
so for a rich enough space $\mathcal{F}$, the procedure captures --and removes from $x$-- the conditional mean $\bar{x}(z)$, i.e. it performs [nonlinear] regression, as in the ``poorest man solution'' of \cite{TT2}.
When the columns of $G$ span a general quadratic function of $y$, not only the conditional mean but also the conditional covariance matrix of $y$ is independent of $z$, and (\ref{X(y)}) implies that the relation between $x$ and $y = T(x, z)$ is linear (with $z$-dependent coefficients), as in the ``poor man solution'' of \cite{TT2}. 
We therefore use for the first stage of our algorithm a $G(y)$ whose columns span all quadratic functions of $y$, to capture the conditional mean and covariance matrix of $x$, leaving a more detailed characterization of $\rho(x|z)$ to subsequent stages.}

One could extend the choices above and fill all columns of $G$ with externally provided functions, such as Hermite polynomials of a given degree. Similarly, we could use as columns of $F$ polynomials when $z \in R^d$ , trigonometric functions when $z$ is a periodic variable and indicator functions when
$z$ can only adopt a finite set of categorical values. Yet it is generally preferable to use a less parametric approach and let the data dictate the form of the functions to use. 
For our experiments, we have used a simple class of data-adapted spaces described in the appendix, where the columns of $F$ and $G$ are given by asymmetric kernel-like functions with column dependent center and bandwidths, a flexible and economic variation of reproducing kernel Hilbert spaces:
$$ F^j(z) = K^z\left(z, z^c_j\right), \quad G^j(z) = K^y\left(y, y^c_j\right). $$
Even though the $\{y_i\}$ evolve, the centers $\{y^c_j\}$ are fixed throughout each stage of the procedure, so as to have a fixed functional space $\mathcal{G}$. Their cardinality does not need to match that of the $\{y_i\}$, it is typically much smaller.

\esteban{Since selecting functional spaces for each stage involves too many hyper-parameters to fit manually, one may consider an adaptive methodology that adjusts them automatically. We leave such development though to further work, since it would extend this article too far beyond its current scope.
A further consideration in this regard is that, ideally, the spaces $\mathcal F$ and $\mathcal G$ should be selected jointly. For instance, regions in $z$ space where $y$ changes rapidly need to be better resolved than those where $\rho(y|z)$ is comparatively flat.
}

\subsection{Successive barycenter problems}\label{sec:FandG}

Since the barycenter problem removes from $x$ any $z$-dependence detectable through the functional spaces $\mathcal{F}$ and $\mathcal{G}$, $y = T(x, z)$ can still depend on $z$ in ways that $\mathcal{F}$ and $\mathcal{G}$ do not capture. For instance, if $\mathcal{G}$ consists only of quadratic function of $y$, just the conditional mean and covariance of $x$ are removed, leaving in $y$ any other $z$-dependent property of $\rho(x|z)$, such as higher moments or the distribution's modality. Similarly, if $\mathcal{F}$ includes only functions of a subset of the $\{z^l\}$, $y$ may still depend on the remaining ones, if $\mathcal{F}$ includes only functions of the individual $\{z^l\}$, any non-additive dependence of $x$ on the $\{z_l\}$ will remain in $y$, and if the bandwidths of the functions in $\mathcal{F}$ are large, only long-scale trends are removed, leaving small-scale signals unresolved.

\esteban{
As pointed out in the prior subsection, one can address this by updating the functional spaces at each stage of the procedure. One must make sure though that each functional space includes that of the prior stage, for else $y$ may return back toward $x$, as some of the constraints that had pushed it toward its current location will not be active anymore.
An alternative that does not suffer from this drawback, is to run not successive stages but successive runs of the OTBP, where the variables playing the role of $x$ in the $l$'th run are the final $y$ from the previous run. Since the cost function of each run does not refer to the original $x$, there is no component of the objective function pulling $y$ back to $x$, undoing part of the progress of the prrvious run.
}
The final $y$ resulting from this multi-run procedure are not samples of the barycenter of the original $x$, since the composition of optimal maps is not necessarily optimal. However, since we know how to invert each of the maps, we can still simulate and estimate by composition $\rho(x|z_*)$ for any target $z_*$.
This procedure resembles boosting \cite{freund1996experiments}, in which multiple models are trained sequentially  so that each new problem removes further variability from the barycenter of the prior one.

\subsection{Complexity}
\label{sec:complexity}
One major advantage of the new methodology is its efficiency, which makes it applicable to large data sets. This efficiently derives from formulating the independence conditions between $y$ and $z$ in terms of the uncorrelation between test functions and subsequently relaxing it to the vanishing of the singular values of a matrix $A(y)$ whose rank does not depend on the sample size.   
Previous methods  \cite{tabak2022distributional} used kernels where every data point acted as a center, yielding at least O$(n^2)$ time complexity. By contrast, the new algorithm's time complexity scales bilinearly with the number of samples ($n$) and the dimension of the data ($d_x$). 
\andrew{It scales also with the user-defined parameter $J_{\max}$, the maximum number of features used at any step, which
depends on the structure of the problem rather than its size. 
}

The algorithm's operations can be broken into those performed only once per run, those performed once per stage, and those performed at each descent step. Even though it is only the third category that determines the time complexity of the algorithm in practice, we analyze all three for completeness.

Since the factors $z$ do not evolve through a run, the orthogonal matrix $Q_z$ is computed only once, \esteban{so its effect on the total computation cost is small except for high-dimensional $z$.}  
When using kernels, calculating $Q_z$ requires k-means clustering to determine the centers for $z$, which with a fixed maximum number of iterations requires O($n*d_z*m_z$) operations. Evaluating the kernel function also requires O($n*d_z*m_z$) operations. A standard singular value decomposition of the matrix $F \in \R^{ n\times n_z}$ requires O($n*n_z^2$) steps where $n_z$ is a user's provided input. Therefore, the number of operations performed only once scales as O($n*d_z*m_z$). For very large data-sets with high-dimensional factors $z$, this number can be further reduced by adopting state-of-the-art methodologies for finding the first few principal components of large matrices \cite{rokhlin2010randomized,artac2002incremental}. 

To recall, a new stage is started when the average squared distance between the current values of the $\{y_i\}$ and their values $\{y_i^0\}$ at the start of the current stage is larger than a prescribed fraction of the variance of the latter. At the beginning of each stage, the orthogonal matrix $Q_y$ needs to be computed. The same scaling arguments apply here as for the calculations of $Q_z$, yielding a total of O($n*d_y*m_y$) operations per stage.

The main loop iteration requires calculating derivatives of the cost function and of the penalty function. The complexity of the former is O($n*d_y*m_y$). The latter requires calculating the gradient of the matrix $G(y)$ which, when using kernels, involves calculating the kernel in $y$ and its derivatives, with complexity O($n*d_y*m_y$) if the number of kernel centers is fixed. This is followed by a matrix multiplication which is O($n$). So overall each iteration performs O($n*d_y*m_y$) operations, which decouple among the sample points, making them trivially parallelizable.

Not captured by the complexity analysis above is the number of iterations required for convergence. 
Additionally, in practice one may adopt larger values of $n_y$ and $n_z$ for problems with more complex dependence between $x$ and $z$. Yet, for a fixed problem, the number of iterations should not depend on the number of data points $n$, an observation confirmed in our numerical experiments.
We verify the algorithm's complexity by plotting the time of the pre-calculations, stage calculations, average descent iteration and total time as $n$, $d_x$ and $d_z$ vary. In each case, each factor $z^l$ is a normal random variable with mean $0$ and variance $0.25$ and $x$ is drawn from the $z$-dependent isotropic gaussian
$$x\sim \mathcal{N}\Bigg(\cos\Big(2\pi\sum_{l=1}^{d_z}z^l\Big)\vec{1}_{d_x}, 0.05\Big[\big[\sin\Big(0.1*\big(\sum_{l=1}^{d_z}z^l+0.2\big)^2\Big)+0.25\big]^{-1}\Big]I_{d_x}\Bigg).$$

Figure 1 displays the data and barycenter for $d_x = d_z = 1$ and figure 2 shows the running times for various values of $n$, $d_x$, and $d_z$. Each data point displayed is the median across 10 trials of the mean time spent in each portion of the algorithm. When the dependence on $n$ is being considered, $d_x$ and $d_z$ are both kept at 1. When either $d_x$ or $d_z$ are being varied, the other is kept at 1 and $n$ is kept at 2000. The first two rows use kernel based functions for both $F$ and $G$ while the third row uses only linear and quadratic functions of $y$ for $G$, which has the same complexity but with lower constants and thus is much faster. 
These experiments confirm that the complexity of all stages grow linearly with the number of samples and that the descent steps contribute most heavily to the total run time. They also show that, as predicted, the dimension of $x$ increases linearly the time complexity of the main descent steps and the calculations per stage, and the dimension of $z$ increases linearly the complexity of the pre-calculations.

\begin{figure}[h!]
       \centering
       \includegraphics[width = 0.3\linewidth]{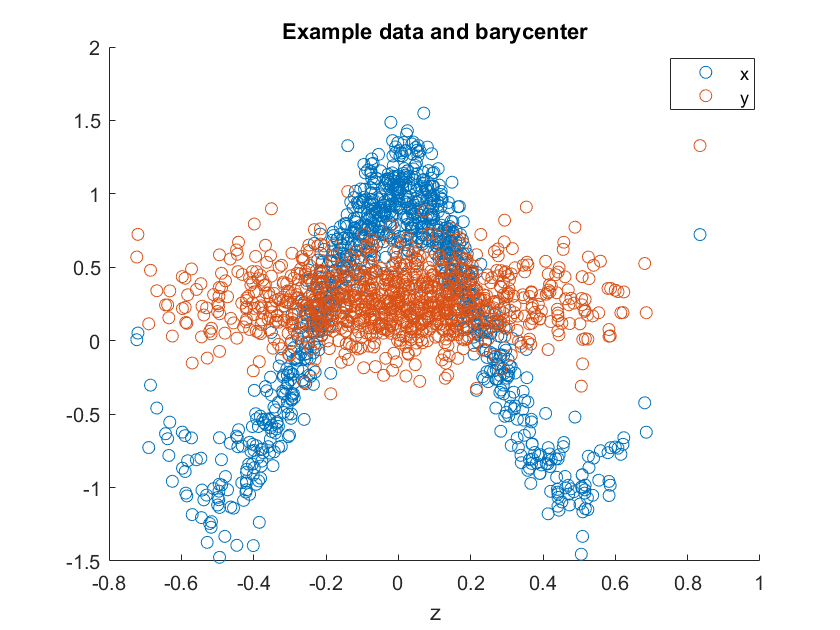}
       \caption{Data points $\{x_i\}$ and corresponding barycenter samples $\{y_i\}$ with $d_x = d_z = 1$. }
       \label{fig::complexity_data}
       \vskip-0.5cm
\end{figure}

\begin{figure}[h!]
        \centering
        \includegraphics[width=0.99\linewidth]{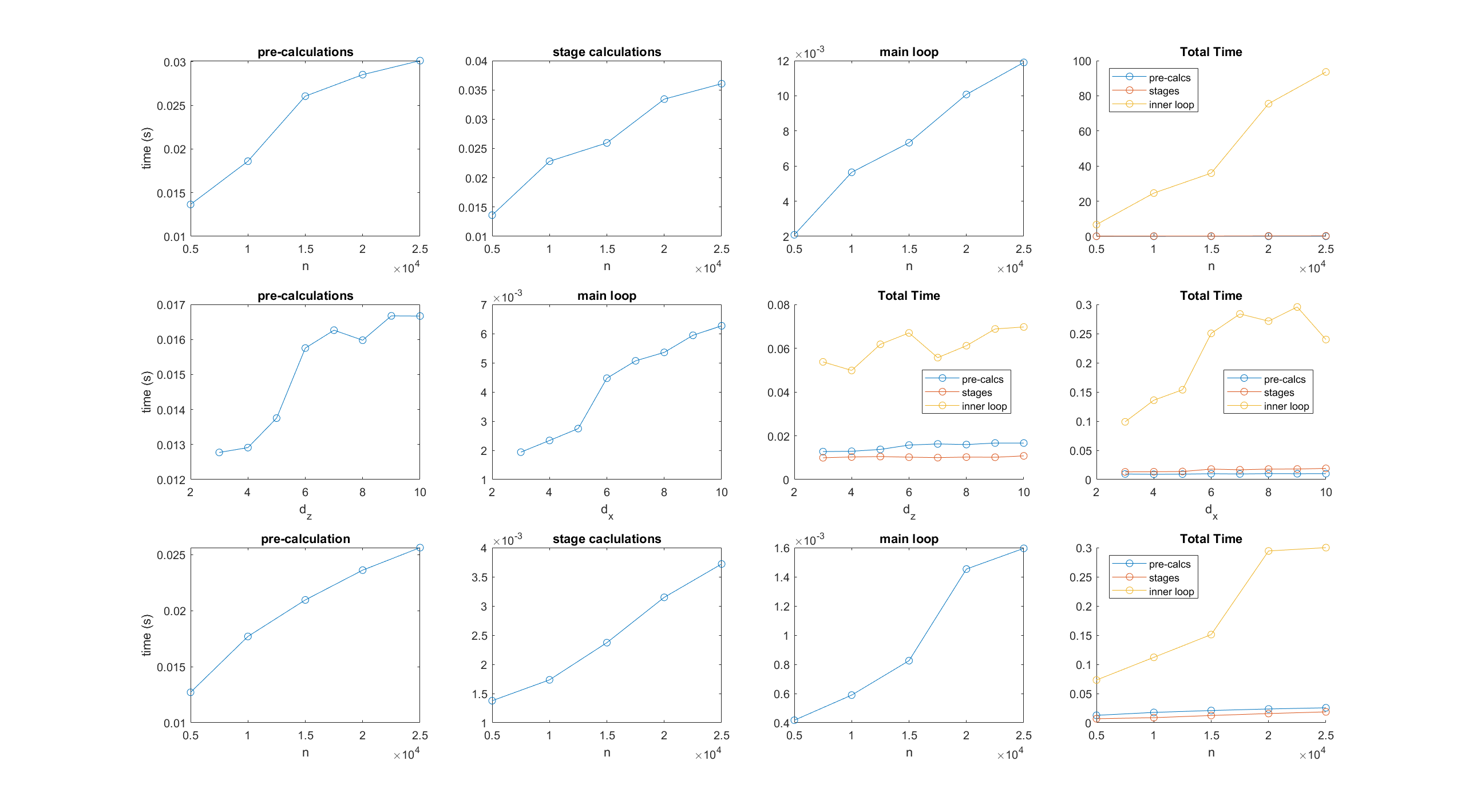}\quad
        \caption{Median time dependency on number and dimension of data. The top and bottom rows display the run-time dependence of the pre-calculations, calculations per stage, calculations per time-step and total run-time, on the number of samples $n$, with $d_x = d_z = 1$, the top one using kernels in $y$, the bottom only linear and quadratic functions.  The middle row displays the relevant dependences on $d_x$ and $d_z$, with the other fixed at $1$ and $n$ at $2000$.}%
        \label{fig::Complexity}%
\end{figure}
\vskip-1cm

\section{The optimal transport problem}\label{sec:OTMonge}

The methodology developed for the OTBP can be easily adapted to solve regular [Monge] OTP.
The OTP is simpler than the OTBP, as it involves only two distributions, a source $\rho_0$ and a target $\rho_1$:
\begin{equation}
 \min_{w = Q(x)} E_{\rho_0} \left[c\left(x, w\right)\right], \quad Q\# \rho_0 = \rho_1. 
 \label{OT}
\end{equation}
 For the purpose of relating them to the barycenter problem, we introduce a binary covariate $z \in \{0, 1\}$, and think of the source and target distributions as instances of the single conditional distribution $\rho(x|z)$:
$ \rho_0(x) = \rho(x|z=0), \quad \rho_1(x) = \rho(x|z=1)$.
Then the map $y=T(x, z)$ that solves the barycenter problem for $\rho(x|z)$, automatically provides the solution to (\ref{OT}) through
$ Q(x) = T^{-1}\left(T(x, 0), 1\right)$,
a standard result in interpolation displacement \cite{santambrogio2015optimal}.

In the data-driven case, we have $n_0$ samples $\{x_i^0\}$ from $\rho_0$ and $n_1$ samples $\{x_j^1\}$ from $\rho_1$, for a total of $n = n_0 + n_1$ pairs $\{x_i, z_i\}$ from $\pi(x, z)$. Since our methodology provides all the values $y_i = T(x_i^0, 0)$ and an explicit formula for $T^{-1}\left(y_i, 1\right)$, we have direct access to all $Q(x_i^0)$ and, mutatis mutandis, we have also access to its inverse, $Q^{-1}(x_j^1)$. Under  the canonical cost, the corresponding formula simplifies to
$$ Q(x_i^0) = T\left(x_i^0, 0\right) - \frac{n_0}{n_1} \left(x_i^0 - T\left(x_i^0, 0\right) \right) , $$
as follows from the fact that every point in the barycenter is the weighted geometrical $c$-barycenter of its pre-images \cite{kuang2019sample}. Then, with only two distributions, a point at the barycenter and one of its pre-images suffice to find the other.

The fact that there is only one, binary covariate $z$ simplifies our methodology considerably, since except for an arbitrary sign, there is  only one function $f(z)$ with zero mean and norm one:
$ f(z) \propto \begin{cases} \ \frac{1}{n_0} & \hbox{for } z=0 \cr -\frac{1}{n_1} 
& \hbox{for } z=1 \end{cases} $.
Then the barycenter problem reduces to
\begin{equation}
 \min_{y_i} L = \sum_{i=1}^n c\left(x_i, y_i\right) + \lambda \left\|f' Q_y(y) \right\|^2,
 \quad f_i = f(z_i),
 \label{barycenter_two}
\end{equation}
where we have used the fact that, since $f$ is fixed, the maximizing vector $b$ replacing the right principal component of $A$ is proportional to ${Q_y}' f$:
$$ \arg \max_{\|b\| = 1} f' Q_y b = \frac{Q_y' f}{\|Q_y' f\|} \Rightarrow \sigma \defeq \max_{\|b\| = 1} f' Q_y b = \|f' Q_y\|. $$
Other than the simplifying facts that we do not need to update $a$ and $b$ and the matrix $Q_z$ consists of a single column, the procedure to solve (\ref{barycenter_two}) follows the same steps as the one for the full barycenter problem (\ref{L_multiK_2}).

We can bypass the barycenter $\mu$ in the procedure above, finding a map $Q(x) = T(x, 0)$ that pushes forward $\rho_0$ to $\rho_1$ directly, by minimizing $L$ in (\ref{barycenter_two}) only over the $y_i$ with corresponding $z_i = 0$, i.e. over $T(x, 0)$, leaving the remaining $y_i$ fixed at $x$, i.e. setting $T(x, 1) = x$. This enforces the condition that $T(x, 0)\# \rho_0 = \rho_1$, since the final 
$ y_i^* = \begin{cases} T\left(x_i, 0\right)& \hbox{for } z_i = 0 \cr  
x_i & \hbox{for } z_i = 1 \end{cases}$
must be independent of $z$. 
This procedure, while lacking the symmetry of the prior one with respect to $\rho_{0,1}$, is more straightforward, and is particularly well-suited for density estimation.

\section{Conditional density estimation}\label{sec:CDE}

Our methodology simulates $\rho(x|z)$, by producing $n$ samples $\{x^*_i\}$ from $\rho(x|z_*)$ for any target $z_*$. Simulation is at the core of many applications, but others, such as Bayesian inference, require the evaluation of $\rho(x|z)$ for arbitrary values of $x$ and $z$. The fact that typically there is none or at most one observation available for any target value $z$ makes estimating $\rho(x|z)$ directly from the data $\{x_i, z_i\}$ challenging. A slight extension of our procedure produces such conditional density estimation.

Regular --as opposed to conditional-- density estimation can be obtained through the OTP as follows. Given $n$ samples $\{x_i\}$ drawn from the unknown distribution $\rho(x)$ that we seek to estimate, select a target distribution $\mu(y)$ that one can easily both evaluate and sample, such as a Gaussian, and find the optimal map $Q(x)$ pushing forward $\rho$ to $\mu$. Then
$$ \rho(x) = \left|\det\left(\nabla_x Q\right) \right| \  \mu(Q(x)) , \quad  \rho(X(y)) = \frac{1}{\left|\det\left( \nabla_y X \right) \right|} \  \mu(y) , \quad
 X = Q^{-1}, $$
so for any $x$,
$$ \rho(x) = \frac{1}{\left|\det\left( \nabla_y X(y)\big|_{y = Q(x)}\right)\right|} \  \mu(Q(x)) . $$
In our procedure, $Q(x) = T(x, 0)$, and $\nabla_y X(y)$ is known from (\ref{Xy(y)}). 

If the density $\rho(x)$ is sought for values of $x$ different from the $\{x_i\}$, one can carry these values through the procedure to their final $y = Q(x)$ as passive tracers that do not affect $L$ in (\ref{barycenter_two}). Alternatively, one can solve the reciprocal OTP from $\mu$ to $\rho$, and then write
$ \rho(x) = \left|\det\left( \nabla_x Y(x)\right)\right| \  \mu(Y(x))$ 
for any value of $x$ sought.

This density estimation procedure requires selecting a target measure $\mu = \rho_1$. We can either adopt a fixed target, such as a standard Gaussian, or adapt it to the data, using for instance a Gaussian with the same mean and covariance matrix as the data or a Gaussian mixture fitted to the data through Expectation Maximization. The advantage of such more tailored approaches is that the corresponding OTP becomes easier, since even the trivial map $Q(x) = x$ provides a regular parametric density estimation. 

\medskip

We can apply this procedure to conditional density estimation, i.e. estimate $\rho(x|z)$ from $n$ samples $\{x_i, z_i\}$ in at least two distinct ways:
\begin{enumerate}
\item Obtain $n$ samples $\{x^*_i\}$ from $\rho(x | z_*)$ and apply density estimation to these directly.

\item Estimate the density $\mu(y)$ of the barycenter, and then compute
$$ \rho\left(X(y, z) | z\right) = \frac{1}{\left| \nabla_y X(y, z) \right|} \  \mu(y), $$
with $X(y, z)$ given by (\ref{X(y)}) and $\nabla_y X(y, z)$ by (\ref{Xy(y)}). 

\end{enumerate}
One would choose the first approach when seeking $\rho(x|z)$ for many values of $x$ and only a handful of values of $z$, and the second when exploring the dependence of the conditional density on $z$, as in Bayesian inference.

\section{Numerical examples}\label{sec:NumEx}

We illustrate the methodology through numerical examples, using both synthetic and real data.

\subsection{Statistical Validation of Conditional Samples}
First, in order to demonstrate that the map inversion described in section \ref{sec:mapinv} generates samples that are statistically similar to samples drawn from the true conditional distribution underlying the data, we generate $n$ data points from the model in section \ref{sec:complexity}, with $d_x = d_z = 1$. After solving the barycenter problem, we use the map inversion formula \ref{X(y)} to simulate $n$ samples of $\rho(x | z=0.2)$. Since the true $\rho$ is a normal distribution with known mean and variance, we use a z-test to test the quality of the samples. We use $100$ randomly sub-sampled points from the $n$ simulated samples, so any change in the p-values is explained by the quality of the sampling rather than the number of samples. For robustness, we do this sampling and z-test 100 times and plot the mean p-value with error bars. We expect the true p-values to be uniformly distributed on [0,1], so we plot the expected values as well as the results from sampling using Nadaraya-Watson kernel regression with Silverman's rule of thumb for the bandwidths.

\begin{figure}[h!]
\centering
\includegraphics[width = 0.7\linewidth]{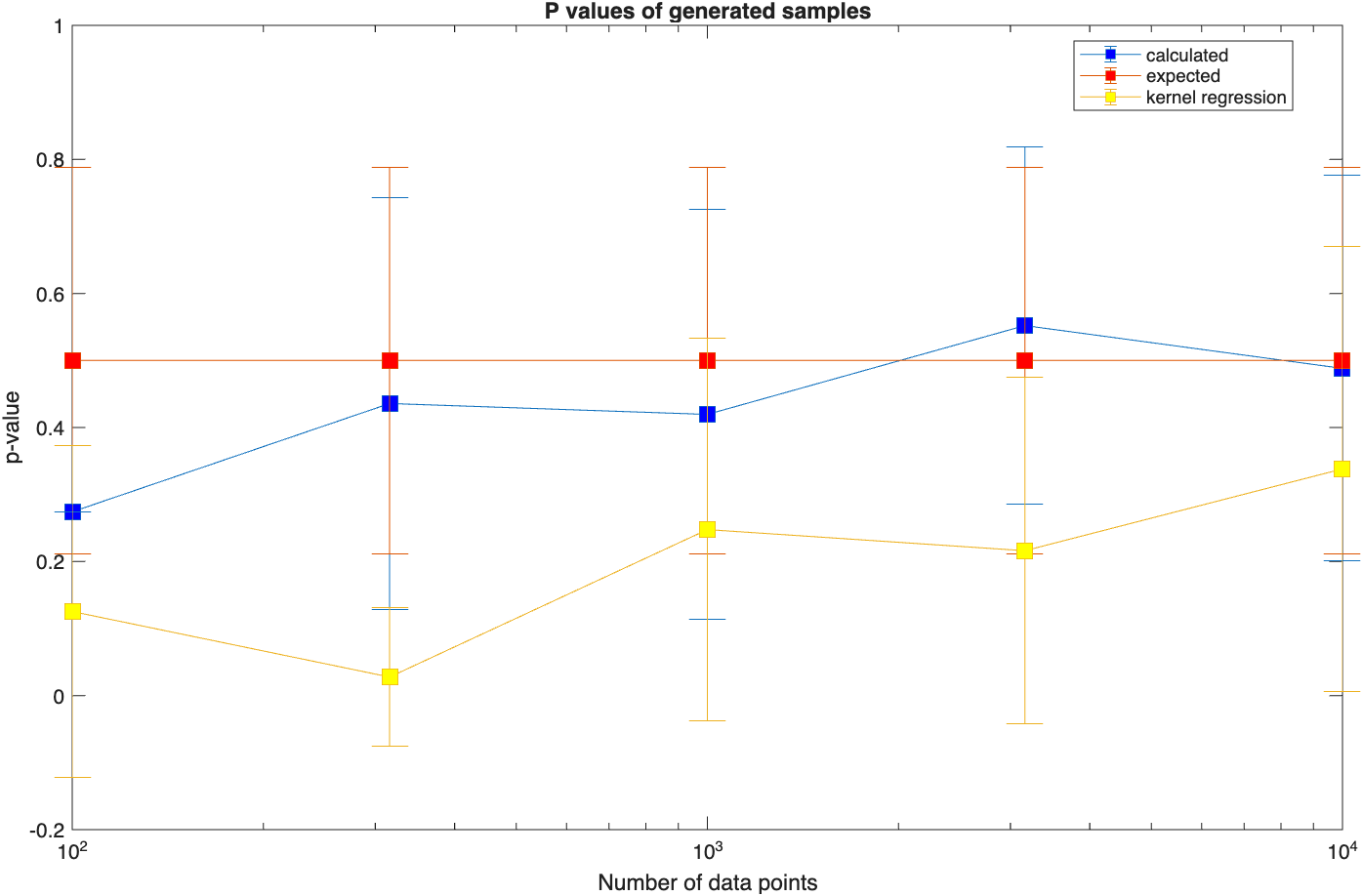}
\caption{distribution of p-values of samples generated with map inversion compared to their expectation and kernel regression}
\label{fig::pvalues}
\vskip-0.5cm
\end{figure}

The figure demonstrates that the simulation not only improves with more data as expected, but that it also does not require many data points to generate good sample data (at least for this specific problem). It outperforms kernel regression, a standard method of conditional density estimation approximation.

\subsection{A Gaussian distribution with $z$-dependent mean and variance}

\begin{figure}[htbp]
    \centering
        \centering
        \includegraphics[width=0.7\linewidth]{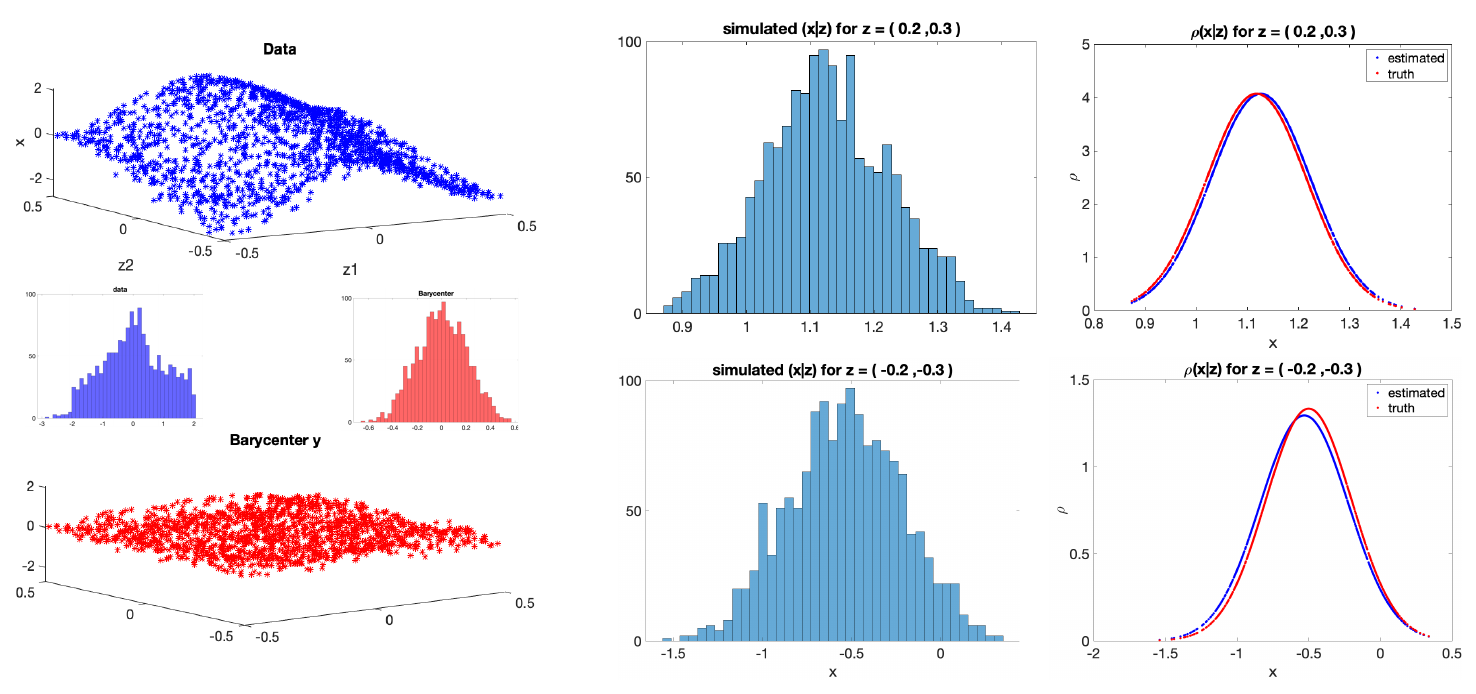} %
        \caption{One-dimensional Gaussian. The leftmost column displays the samples $\{x_i\}$ in the top row and the $\{y_i\}$ in the barycenter in the bottom row, with corresponding histograms showing that the final barycenter converges to a Gaussian. The middle and rightmost columns show the simulated samples $\{x_i^*\}$ and the estimated versus the true density $\rho(x | z_*)$, in blue and red respectively, for $z_* = (0.2,0.3)$ and $z_* = (-0.2,-0.3)$.}
        \label{fig::Gauss1d}
\end{figure}

As a second example, we draw 1500 independent samples from the distribution
$$ x \sim N\left(\mu(z), \sigma^2(z)\right), \quad \mu(z) = \cos\left(2 \pi z_1\right) + \sin\left(\pi z_2\right),
\quad
\sigma(z) = 0.2 \sqrt{\left(1 - 2 z_1\right)\left(1 - 2 z_2\right)}, $$
with $z = \left(z_1, z_2\right)$ uniformly distributed in the square $ -\frac{1}{2} \le z_{1, 2} \le \frac{1}{2}$,
displayed on the top left panel of figure \ref{fig::Gauss1d}. Since for each value of $z$ the distribution for $x$ is Gaussian, it can be fully captured using the two-dimensional test function space $\mathcal{G}(y)$ spanned by the functions $y$ and $y^2$, while keeping for $\mathcal{F}(z)$ a general adaptive space based on kernels. We display the results of the run in figure \ref{fig::Gauss1d} through the corresponding $\{y_i\}$
and the simulation and estimation of $\rho(x | z_*)$ for two selected values of $z_*$, with the true underlying distribution also drawn for comparison. 

We use two measures to quantify the accuracy of the estimated density $\rho_{est}$, in this and other cases where the true density $\rho_{ex}$ underlying the data is known: the relative entropy or Kullback-Leibler divergence (KL)
$$ KL\left(\rho_{est}, \rho_{ex}\right) = \int \log\left(\frac{\rho_{est}(x)}{ \rho_{ex}(x)}\right) d\rho_{est}(x) \approx \frac{1}{n} \sum_{i=1}^n \log\left(\frac{\rho_{est}(x_i)}{ \rho_{ex}(x_i)}\right),$$
where the $\{x_i\}$ are the $n$ samples generated by the algorithm, and the Mean Square Error (MSE)
$$ MSE = \frac{1}{n} \sum_{i=1}^n \left[\rho_{est}(x_i) - \rho_{ex}(x_i) \right]^2. $$
 Both are reported in table \ref{table_diff}.

\begin{figure}[htbp]
    \centering
        \centering
        \includegraphics[width=0.4\linewidth]{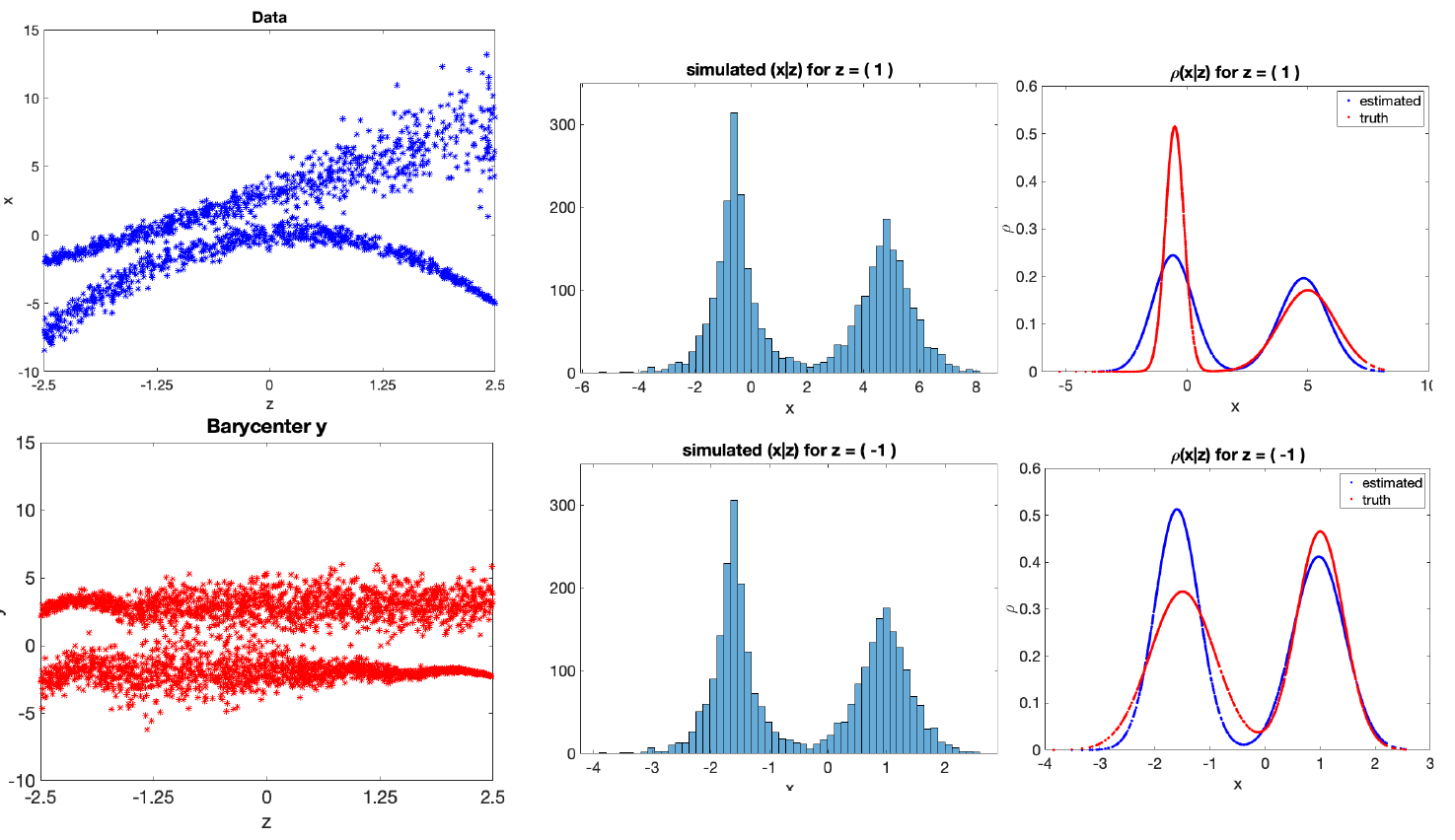} %
    \qquad
        \includegraphics[width=0.4\linewidth]{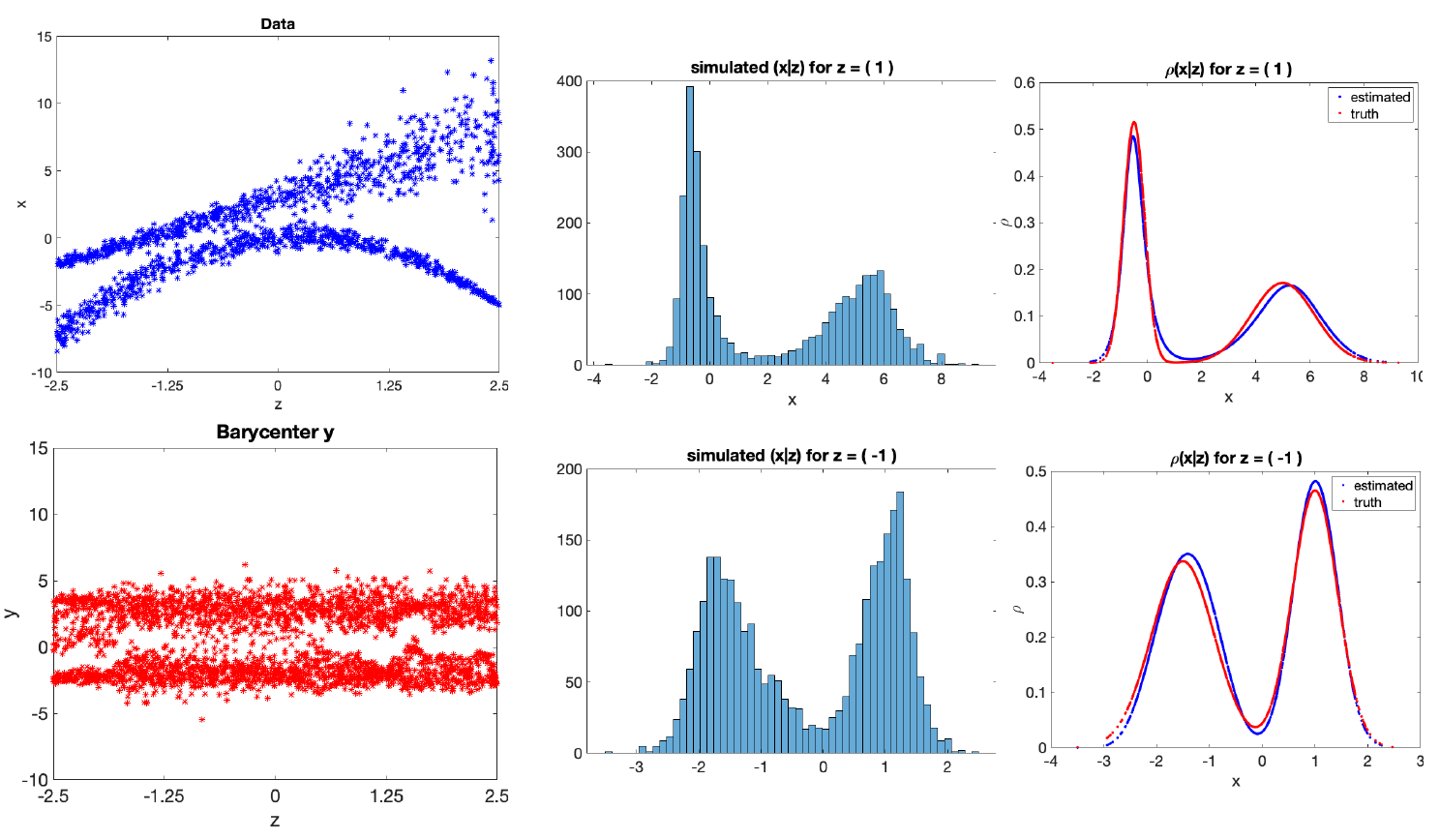} %
        \caption{One-dimensional Gaussian mixture. The left panel displays results when the test function space $\mathcal{G}(y)$ are only linear and quadratic, while the right panels displays results with $\mathcal{G}(y)$ extended to nonlinear kernels. In both panels, the leftmost column displays the data samples in the top row and the barycenter in the bottom row. The middle and rightmost columns show the simulated samples and the estimated versus the true density, in blue and red respectively, for $z_* = 1$ and $z_* = -1$.}
    \label{fig::GMM1d}
\end{figure}

\begin{table}[ht]
\centering
\begin{tabular}{l|l|l|l}
\hline \hline
\toprule
        Example & $z_*$ & MSE & KL-div \\
        \midrule
        \multirow{2}{*}{Gaussian 1D (linear and quadratic $\mathcal{G}(y)$)} & (0.2,0.3)   & 0.0139 & 0.0036 \\
                                     & (-0.2,-0.3)    & 0.0128 & 0.0031 \\
                                     
        \hline \hline %
        
        \multirow{2}{*}{GMM 1D (linear and quadratic $\mathcal{G}(y)$)}   & 1     & 0.0137 & 0.4248  \\
                                     & -1     & 0.0100 & 0.0579 \\
                                                                  
        \hline %
        
        \multirow{2}{*}{GMM 1D (successive kernels $\mathcal{G}(y)$)}   & 1     & $5.9006e^{-4}$ & 0.0776  \\
                                     & -1     & $3.5332e^{-4}$ & 0.0015 \\
                                     
		\hline \hline %
        
		\multirow{2}{*}{GMM 2D (linear and quadratic $\mathcal{G}(y)$)}   & 0.5     & $6.3706e^{-4}$ & 2.6249  \\
                                     & -1     & 0.0012 & 2.5624 \\
                                                                  
        \hline %
                
        \multirow{2}{*}{GMM 2D (successive kernels $\mathcal{G}(y)$)}   & 0.5     & $2.4176e^{-4}$ & 1.7166
  \\
                                     & -1     & $7.2297e^{-4}$ & 1.4013 \\
        
        \bottomrule
 \hline 
\end{tabular}
\caption{Metrics (mean squared error and KL divergence) measuring the difference between estimated densities and true densities, with each row entry corresponding to a $z_*$ value for every set of numerical experiments.}
\label{table_diff}
\end{table}

\subsection{Two $z$-dependent Gaussian mixtures}

In order to consider non-Gaussian examples --more generally, examples where the dependence of $x$ on $z$ cannot be reduced to a $z$-dependent linear transformation-- we perform first a one-dimensional experiment, drawing 1500 samples from the $z$-dependent Gaussian mixture displayed on the top left panel of figure \ref{fig::GMM1d}:
$$ x \sim \sum_{k=1}^2 \gamma_k\ N\left(\mu_k(z), \sigma_k^2(z)\right), \quad \gamma_1 = \gamma_2 = \frac{1}{2}, $$
$$ \mu_1(z) = 3 + 2 z, \quad \sigma_1^2(z) = \frac{1}{2} e^{z} , \quad \mu_2(z) = \frac{z}{2} - z^2, \quad \sigma_2^2(z) = 0.25 - 0.1 z, \quad z \sim \text{U}([-2.5, 2.5]) . $$
The results of using only linear and quadratic functions for $\mathcal{G}(y)$, displayed in the left panel of figure \ref{fig::GMM1d}, show how a linear transformation cannot remove the dependence of $y$ on $z$, and so the corresponding simulation and estimation of $\rho(x|z)$ are far from accurate.
 Then, we performed four successive barycenter problems, the first with just linear and quadratic $\mathcal{G}(y)$, the rest with adaptive kernels, with the bandwidths of the kernels for both $\mathcal{G}(y)$ and $\mathcal{F}(z)$ adopted smaller for each successive run. The results, displayed in the right panel of figure \ref{fig::GMM1d}, show how the simulated samples and conditional density estimation recover the original $z$-dependent Gaussian mixture much more robustly. The metrics in table \ref{table_diff} confirms that a more complex test function space decreases the estimation error.

\begin{figure}[htbp]
    \centering
    \begin{subfigure}[t]{0.48\textwidth} %
        \centering
        \includegraphics[width=0.8\linewidth]{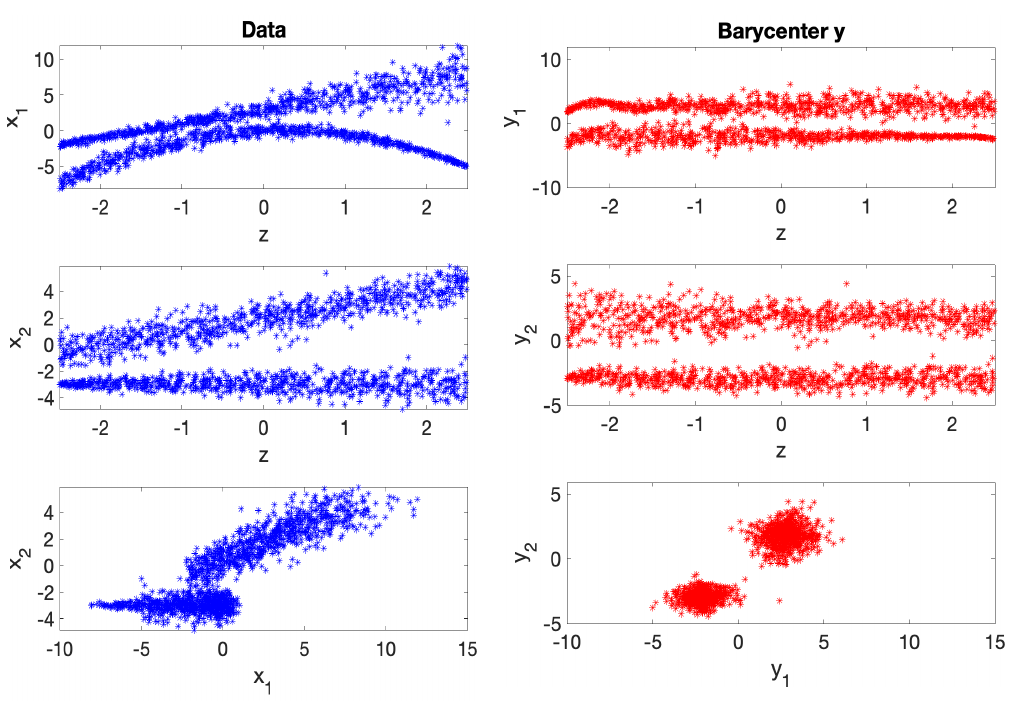} %
        \caption{A two-dimensional Gaussian mixture. The left panel displays the data points $\{x^{1,2}_i\}$ as a function of the corresponding $\{z_i\}$ and the $\{x^1_i, x^2_i\}$, while the right panel shows the equivalent plots for the $\{y^{1,2}_i\}$ in the barycenter.}
    \end{subfigure}
    \hfill %
    \begin{subfigure}[t]{0.48\textwidth} %
        \centering
        \includegraphics[width=\linewidth]{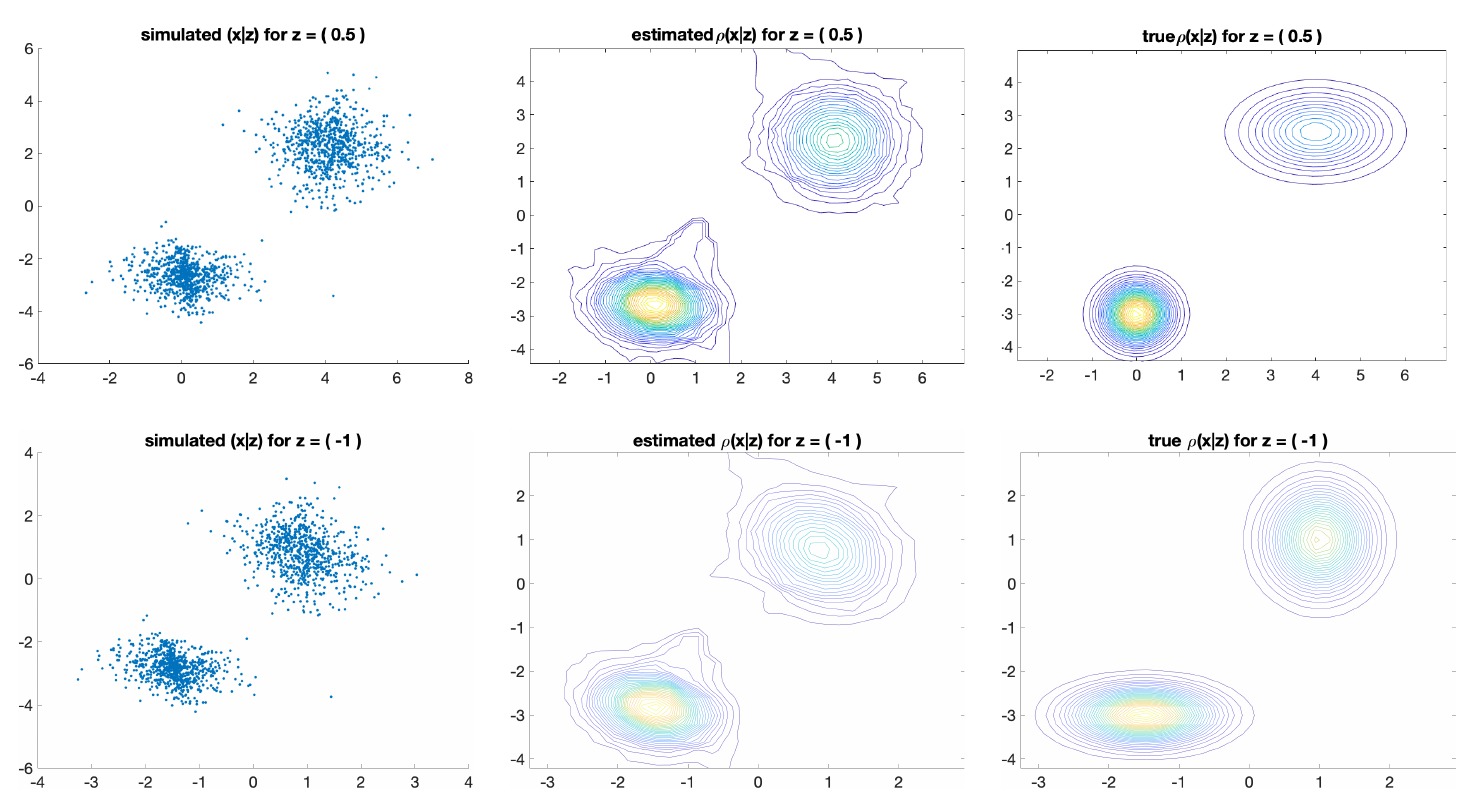} %
        \caption{Two-dimensional Gaussian mixture. The top and bottom rows show the simulated samples, the estimated conditional density, and the true density for $z_* = 0.5$ and $z_* = -1$, respectively.}
    \end{subfigure}
    \caption{Two-dimensional Gaussian mixture.}
    \label{fig::GMM2d_data}
\end{figure}

We extend this example to the two-dimensional Gaussian mixture, visualized in panel (a) of figure \ref{fig::GMM2d_data}:
$$ x \sim \sum_{k=1}^2 \gamma_k\ N\left(\mu_k(z), \Sigma_k(z)\right), \quad \gamma_1 = \gamma_2 = \frac{1}{2} , \quad z \sim  \text{U}([-2.5, 2.5]), $$
$$ \mu_1(z) = \left(\begin{array}{c} 3 + 2 z \\ 2 + z \end{array}\right), \quad \Sigma_1(z) = \left(\begin{array}{cc}\frac{1}{2} e^{z} & 0 \\ 0 & 0.5 \end{array}\right) , \quad \mu_2(z) = \left(\begin{array}{c}\frac{z}{2} - z^2 \\ -3 \end{array}\right), \quad \Sigma_2(z) = \left(\begin{array}{cc}0.25 - 0.1 z & 0 \\ 0 & 0.25 + 0.1 z \end{array}\right) . $$
The results of the procedure are displayed in panel (b) of figure \ref{fig::GMM2d_data}, and estimation error displayed in table \ref{table_diff}.

\subsection{An example of Bayesian inference}

\begin{figure}[h!]
        \centering
        \includegraphics[width=0.18\linewidth]{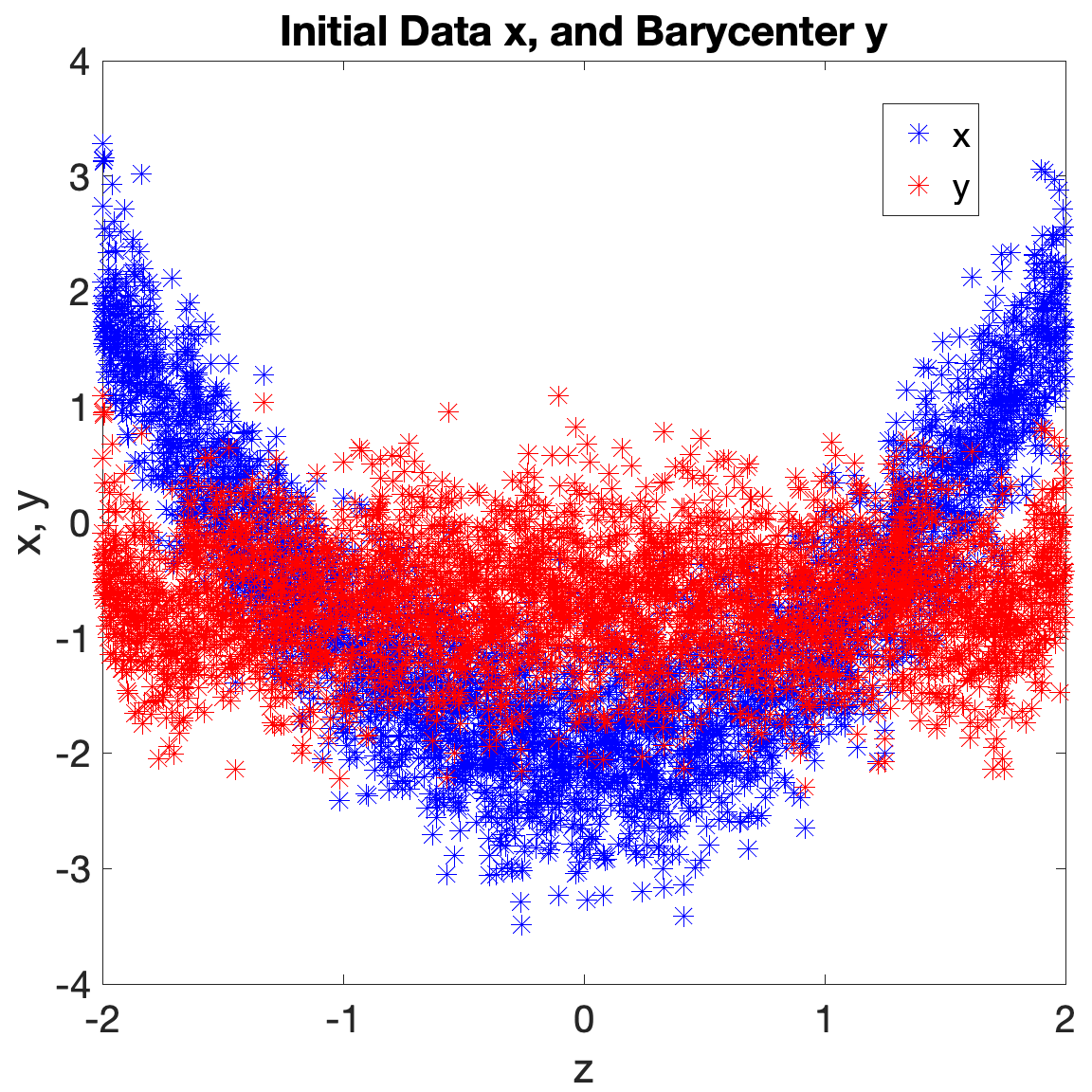}\quad
        \includegraphics[width=0.25\linewidth]{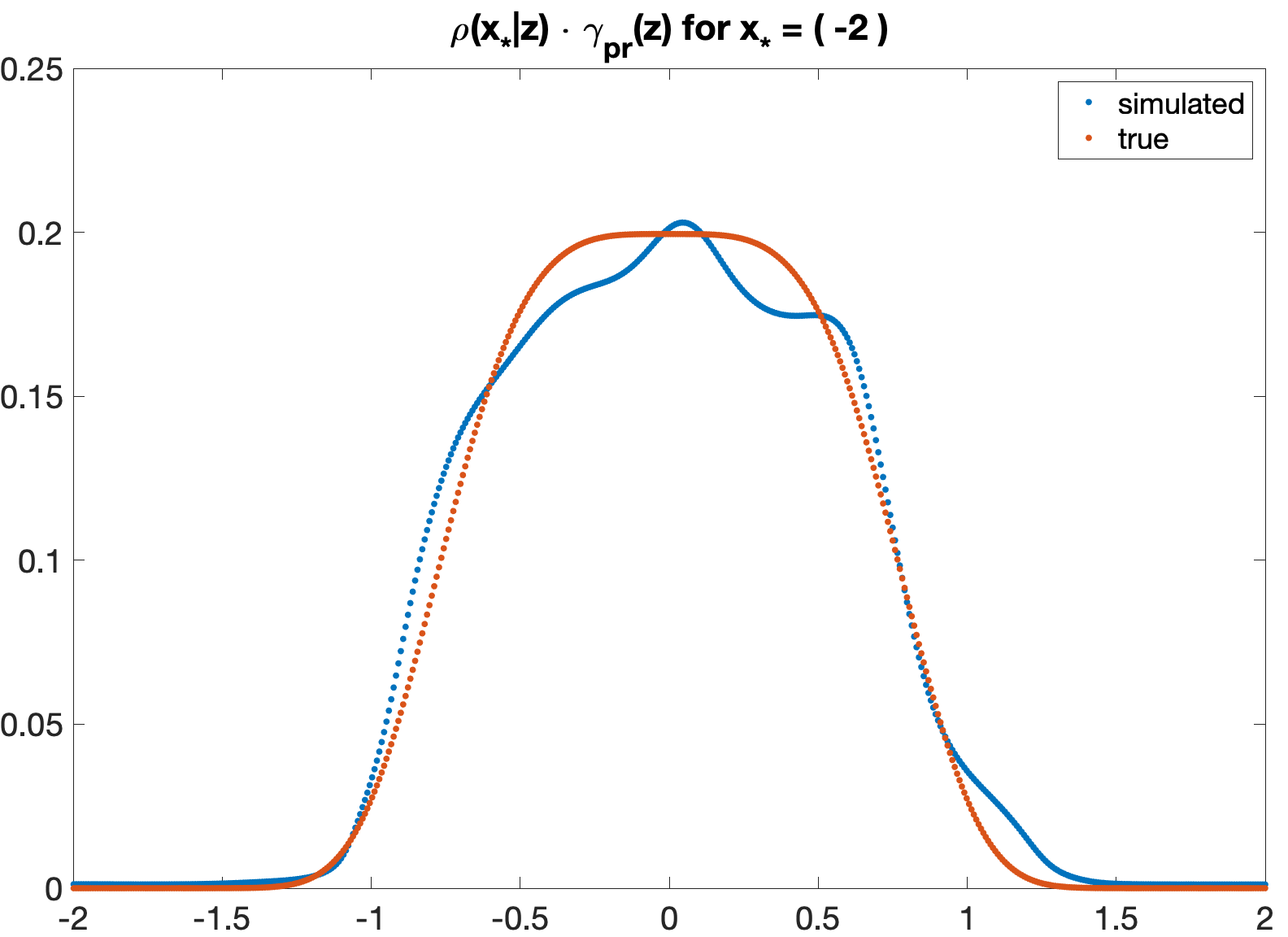}\quad
        \includegraphics[width=0.25\linewidth]{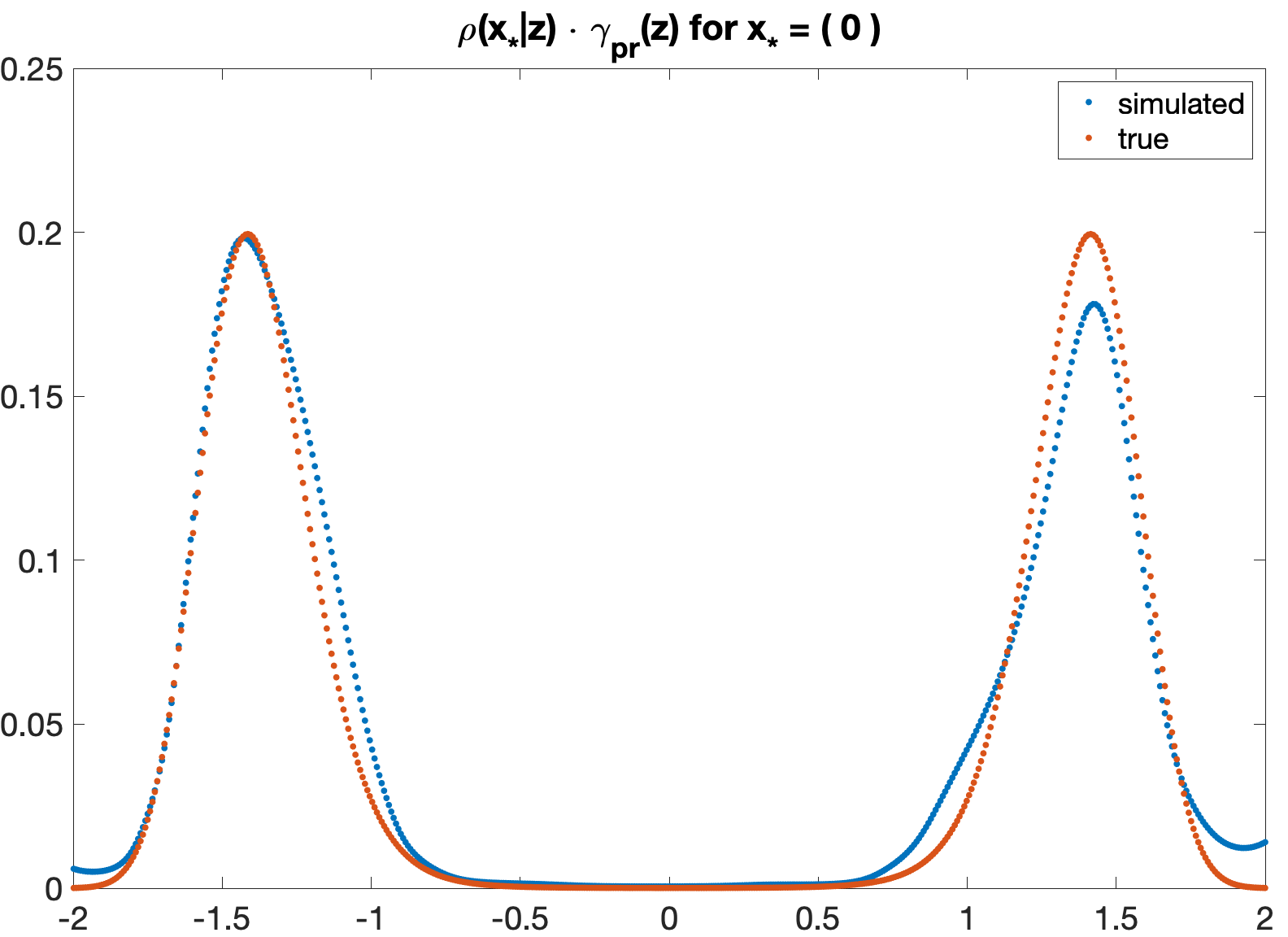}
        \caption{Original data, barycenter and simulated versus true posterior density $\gamma(z | x^*) $ evaluated at $x^*=-2$ and $x^* = 0$.}%
        \label{fig::BI}%
    \end{figure}

This section illustrates model-free Bayesian inference using the OTBP (A different use of push-forward maps for Bayesian inference \cite{ELMOSELHY20127815} pushes forward the prior to the posterior measure.) To demonstrate our approach, we draw samples from the distribution
$
\rho(x|z) = N(z^2-2,\sigma^2)$, $\sigma=0.5$, $z \sim \text{U}([-2,2])$. 
The left panel of figure \ref{fig::BI} displays the data $\{x_i\}$ and the discovered $\{y_i\}$ as functions of the corresponding $\{z_i\}$. 
From these, we can directly infer the distribution of $z$ given an observation $x$:
$$ \gamma(z | x) = \frac{\rho(x |z)}{\rho(x)} \gamma_{pr}(z) \propto \rho(x |z) \cdot \gamma_{pr}(z), $$
where $\rho(x|z)$ is the distribution inferred from the data through the OTBP. 
We have adopted as prior $\gamma_{pr}$ the distribution underlying the observed $\{z_i\}$. The results for two values of $x$ are displayed on the middle and right panels of Figure \ref{fig::BI}, overlapped with the exact answer.
They succeed in capturing the transition from unimodal to bimodal distributions corresponding to the parabolic dependence of the conditional mean of $x$ on $z$.
    
\subsection{Online model estimation in Ornstein–Uhlenbeck processes and Lotka-Volterra pray-predator models with observational noise}

We consider next the online estimation of parameters, a key component of data assimilation. Given successive samples from a time series $x^n$ drawn from some transitional distribution $\rho(x^{n+1} | x^n, w^n, \alpha)$ depending on known and unknown parameters $w^n$ and $\alpha$ respectively, and assuming a prior distribution $\gamma^0(\alpha)$ for the latter, we seek to successively improve on these priors as new observations arrive, using Bayes rule:
$$
\gamma^{n+1}(\alpha) \propto \rho(x^{n+1}|x^{n}, w^n, \alpha) \cdot \gamma^{n}(\alpha),
$$
with the proportionality constant determined by the normalizing condition that $\int d\gamma(\alpha) = 1$. In the conventional setting, the conditional distribution $\rho$ is known except for the parameters $\alpha$. We can extend this framework to situations where $\rho(x^{n+1} | x^n, w^n, \alpha)$ itself is only known from partial observations of time series under different values of $\alpha$ and $w$. In a medical setting for instance, $x$ may represent glucose concentration in the bloodstream, $w^n$ the caloric intake at time $t_n$, and $\alpha$ a patient's parameter that may only be determined after treatment. Having observed in the past a number of patients under different diets and having determined their corresponding parameters $\alpha$, we can use this for the online estimation of $\alpha$ for a patient currently under treatment. 
Notice that, in our notation for the OTBP, $x^{n+1}$ plays the role of $x_{n}$ and the pairs $\{x^{n}, \alpha^n\}$ play the role of $z_n$.

\begin{figure}[h!]
        \centering
        \includegraphics[width=0.7\linewidth]{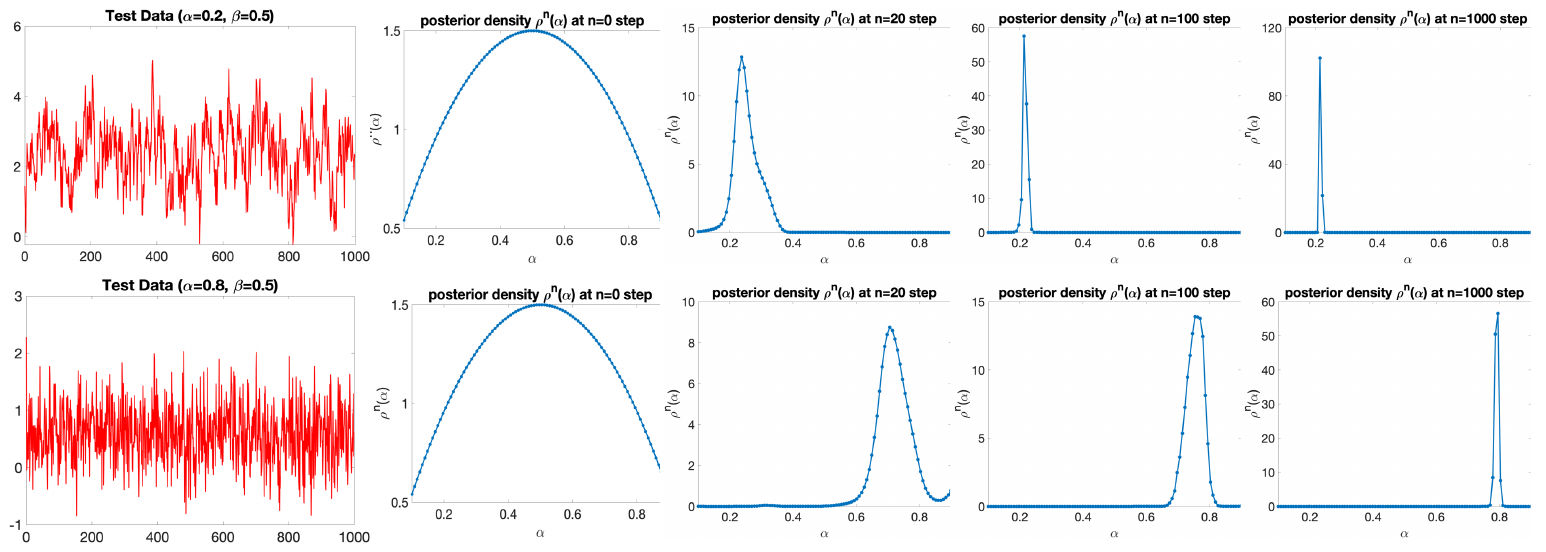}
        \caption{Time-discretized Ornstein–Uhlenbeck process. Testing time series for $\alpha=0.2, 0.8$. Estimated posterior densities as the number of time steps grows for $\alpha=0.2$ in the top row and $\alpha=0.8$ in the bottom row.}%
        \label{fig::OU_DE}%
    \end{figure}

For a first simple example, consider the time-discretized 1D Ornstein–Uhlenbeck process
$$
x^{n+1} = (1-\alpha) x^{n} + \beta + \sigma e^n, \quad e^n \sim N(0, 1) ,
$$
where $\alpha\in (0,1)$ is an unknown model parameter and $\beta = \sigma = 0.5$ are fixed drift and noise levels. Our goal is to learn the model from a set of training data pairs $(x_{\text{train}}^{n+1};\, x_{\text{train}}^{n}, \alpha_{\text{train}}^n)$ and use the model learned to estimate $\alpha$ online from a testing series, while making increasingly more accurate forecasts. For the training data, we draw
$\alpha_{\text{train}}^n$ from a beta distribution $B(2,2)$ over $(0,1)$, which we also adopt as prior $\gamma^{0}(\alpha)$, and $x_{\text{train}}^{n}$ from the uniform distribution $\text{U}([a, b])$.
We carry out experiments with two different parameter values, $\alpha=0.2, 0.8$, with the corresponding  test data displayed in the leftmost column in Figure \ref{fig::OU_DE}. The following results in Figure \ref{fig::OU_DE} demonstrate that the posterior densities converge to delta functions centered around the corresponding true parameters.

In order to apply the procedure to a more complex scenario, we consider next the Lotka-Volterra predator-prey model
\begin{align*}
\frac{d x_1}{dt} &=  \alpha x_1 -\beta x_1 x_2, \\
\frac{d x_2}{dt} &= -\gamma x_2 +\delta x_1 x_2,
\end{align*}
with data observed at irregular discrete times with non-uniform time intervals $\Delta t^n \sim \text{U}[0.5,1.5]$ and corrupted by noise. The time $\Delta t^n$ here plays the role of the known covariate $w^n$, so the Barycenter problem needs to resolve the data dependence on this additional factor.
The data are generated through the explicit trapezoidal numerical scheme but with a much smaller $\tilde{\Delta t}$, to accurately solve the system of ODEs.  Gaussian noise with amplitude  $\epsilon$ is added after simulating the time series to represent noisy observations.

\begin{figure}[h!]
        \centering
        \includegraphics[width=0.7\linewidth]{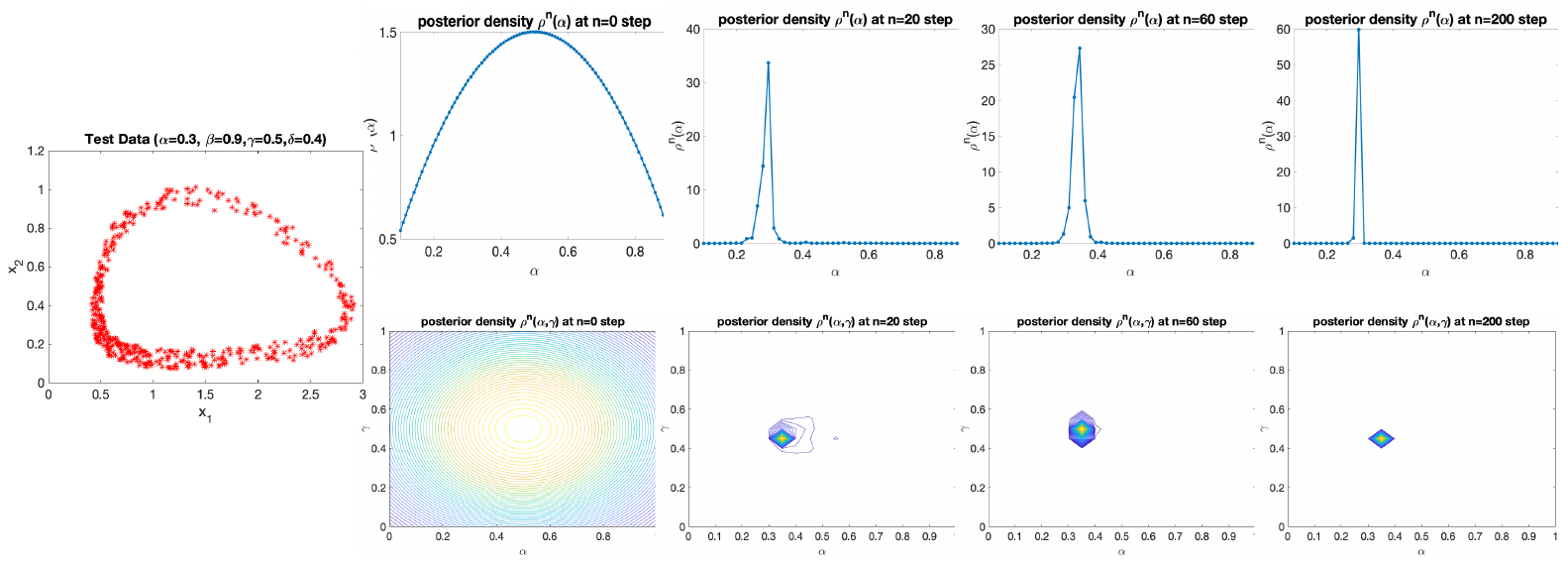}
        \caption{Lotka-Volterra model: Testing time series for $\alpha=0.3, \beta=0.9, \gamma=0.5,\delta=0.4$. Estimated posterior densities as the number of time steps grows for the single unknown $\alpha=0.3$ in the top row, and for two unknowns $\alpha=0.3, \gamma=0.5$ in the bottom row.}%
        \label{fig::LV_DE}%
    \end{figure}

We adopt as test data a simulation with parameters  $\alpha=0.3, \beta=0.9, \gamma=0.5, \delta=0.4$ and $\epsilon=0.1$, yielding the periodic results displayed in the leftmost panel of Figure \ref{fig::LV_DE}.
For a first experiment, we take $\alpha$ as the only unknown parameter.  Again, random training data pairs $(x_{\text{train}}^{n+1};\, x_{\text{train}}^{n}, \alpha_{\text{train}}^n, \Delta t^n)$ are generated, drawing $x_{\text{train}}^{n}$ from a uniform distribution, $\alpha_{\text{train}}^n$ from the beta distribution $B(2,2)$, and deriving the corresponding $x_{\text{train}}^{n+1}$  from the model with additive noise of level $\epsilon$. 
As before, we learn the conditional density $\rho(x^{n+1}|x^{n},\alpha, \Delta t^n)$ by solving the barycenter problem for the training data, with $x_n = x^{n+1}$, $z_n = \{x^n, \Delta t^n\}$. Then we apply Bayes rule online to the testing data, updating at each step 
$$
\rho^{n+1}(\alpha) \propto \rho(x^{n+1}|x^{n},\alpha,\Delta t^n) \cdot \rho^{n}(\alpha).
$$
The results are shown in the top row of Figure \ref{fig::LV_DE}.
We consider next a situation where two parameters, $\alpha$ and $\gamma$, are unknown, so the training data consists of quintuples $(x_{\text{train}}^{n+1};\, x_{\text{train}}^{n}, \alpha_{\text{train}}^n, \gamma_{\text{train}}^n,\Delta t^n)$, and the joint posterior density should be estimated through
$
\rho^{n+1}(\alpha,\gamma) \propto \rho(x^{n+1}|x^{n},\alpha,\gamma,\Delta t^n) \cdot \rho^{n}(\alpha,\gamma)
$, 
with joint Gaussian prior $\rho^{0}(\alpha,\gamma) = N\left(\begin{pmatrix}
0.5 \\ 0.5
\end{pmatrix},\begin{pmatrix}
0.2 & 0 \\ 0 & 0.2
\end{pmatrix}\right)$. The results are displayed in the bottom row of Figure \ref{fig::LV_DE}. We see how in both cases the estimation converges to a delta function centered at the right underlying value of the parameters.

\subsection{Uncovering a hidden signal}

The solution to the barycenter problem helps uncover a hidden signal $w$ that, together with the known factors $z$, fully explain the outcome variable $x$. In order to demonstrate this through examples, rather than simulating a distribution
$ \rho(x | z)$,
we propose a function
$$ x = \phi(z, w), \quad z \sim \gamma(z), \quad w \sim \nu(w) $$
where $w$, playing the role of noise in the distribution, is a hidden cause of variability in $x$.

Recall from Lemma \ref{th:hidden} that the solution $y = T(x, z)$ of the barycenter problem is a proxy for the variable $w$, with $\phi(z, w) \rightarrow X(y, z)$. Moreover, $y$ is related to any ``true'' hidden variable $w$ through a possibly $z$-dependent function
$$ y = Y_z(w), $$
which is invertible if $w$ is identifiable, i.e. if a single value of $x$ cannot originate from a single $z$ and two different values of $w$.

\begin{figure}[h!]
       \centering
       \includegraphics[width = 0.8\linewidth]{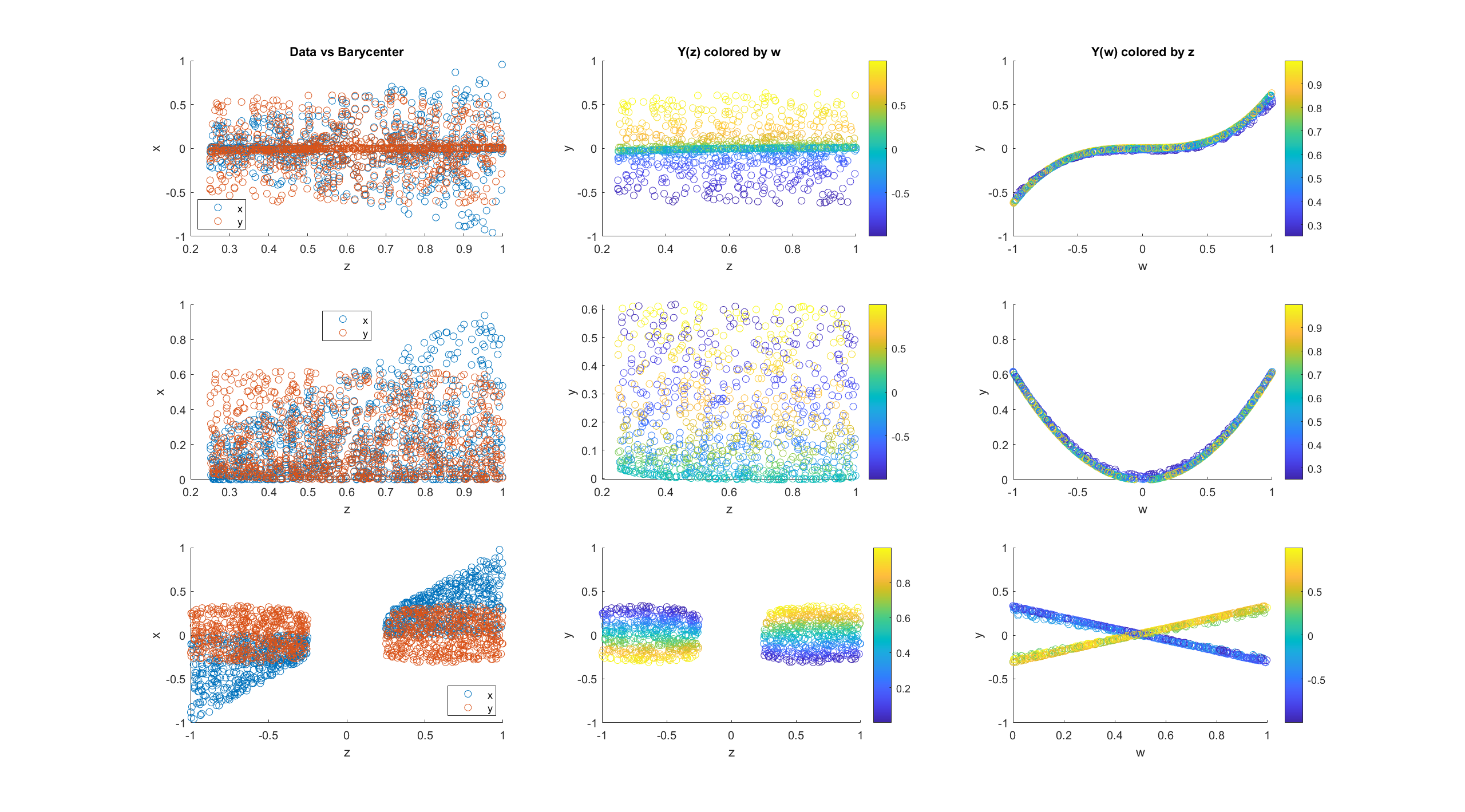}%
       \caption{Three examples with different relations $Y_z(w)$ between $w$ and $y$: one-to-one on the first row, two-to-one on the second and $z$-dependent on the third.}
       \label{fig::HS1}
\end{figure}

Figure \ref{fig::HS1} presents three synthetic examples in order to illustrate the different kind of dependence between $y$ and $w$ typically observed in applications. The panels of each row correspond to the different synthetic examples described below. The left column display the $\{x_i\}$ and corresponding $\{y_i\}$ in terms of the $\{z_i\}$, the middle column displays $y(z)$ again, colored according to the corresponding value of $w$, and the right column displays $y(w)$, colored according to $z$.
In the first row, $z \sim \text{U}[0.25,1]$, $w \sim \text{U}[-1,1]$ and $x = \phi(z, w) = z w^3$ (we exclude values of $z$ near $0$ because $\forall w \ \phi(0, w) =0$, i.e. $\rho(x | 0)$ does not vanish on small sets.) In this example $Y_z(w)=Y(w)$ does not depend on $z$ and $Y(w)$ is invertible. In the second row $x=z w^2$ under the same distributions for $z$ and $w$. We still have that $Y_z(w)$ does not depend on $z$ but now $Y(w)$ is not globally invertible, a reflection of the fact that the sign of $w$ is not identifiable, since $\forall z$ and  $\forall w$ we have that $\phi(z, -w) = \phi(z, w)$. In the third row, $w \sim \text{U}[0,1]$, $z \sim \text{U}[(-1,-0.25)\cup(0.25,1)]$
and $x=z w$,  for which $Y_z(w)$ depends on the sign of $z$.

The analysis of the barycenter underlying the points $\{y_{i}\}$ may at first seem similar to residual analysis, whereby the difference between actual and predicted values is further analyzed to assess model adequacy and improve its predictive power \cite{anscombe1963examination,martin2017fitting}. Both procedures aim to remove variability in the data $x$ attributed to the cofactors $z$, yet while residual analysis only removes $\bar{x}(z)$, the conditional expected value of $x$, the barycenter does this at the level of the full probability distributions $\rho(x|z)$ underlying the data.
Consider the example in the first row of Figure \ref{fig::HS1}, where the data are generated according to $x=zw^3$ with $z$ and $w$ drawn uniformly. The distribution of the residuals (obtained by subtracting the regression of $x$ vs $z$ from the actual value of $x$) would, in this case, be identical to the original distribution underlying the $x$, providing no new useful information. By contrast, the  points $y_{i}$ in the barycenter represent the full variability of $x$ not explainable by $z$. If the $y$'s can be related to known factors $w$, this can be used to improve the model for $x$, for instance by regressing $y$ against these factors and then using both $z$ and the reconstructed $y$ to predict $x$. Better still, instead of regression, the barycenter problem can be used once again to simulate $\rho(y|w)$.

\begin{figure}[h!]
       \centering
       \includegraphics[width = 0.3\linewidth]{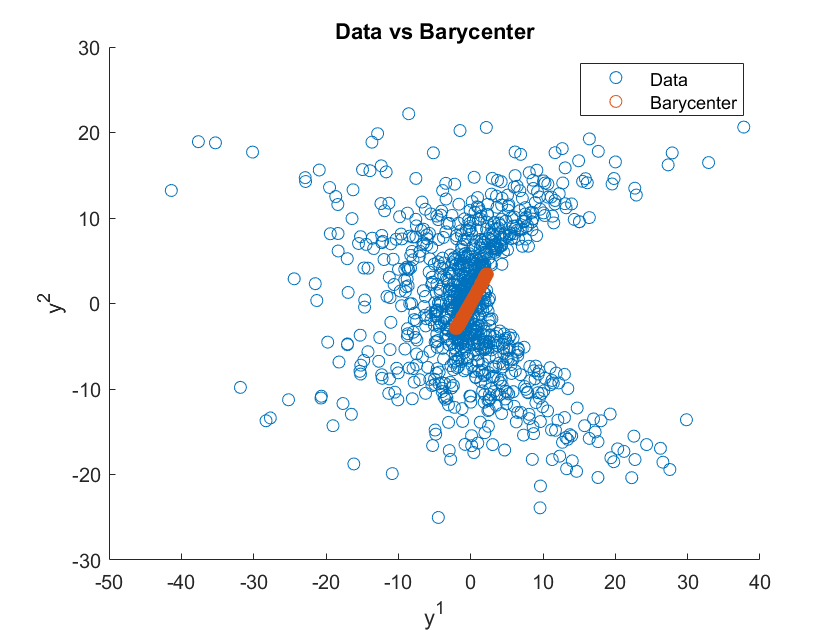} \quad%
       \includegraphics[width = 0.3\linewidth]{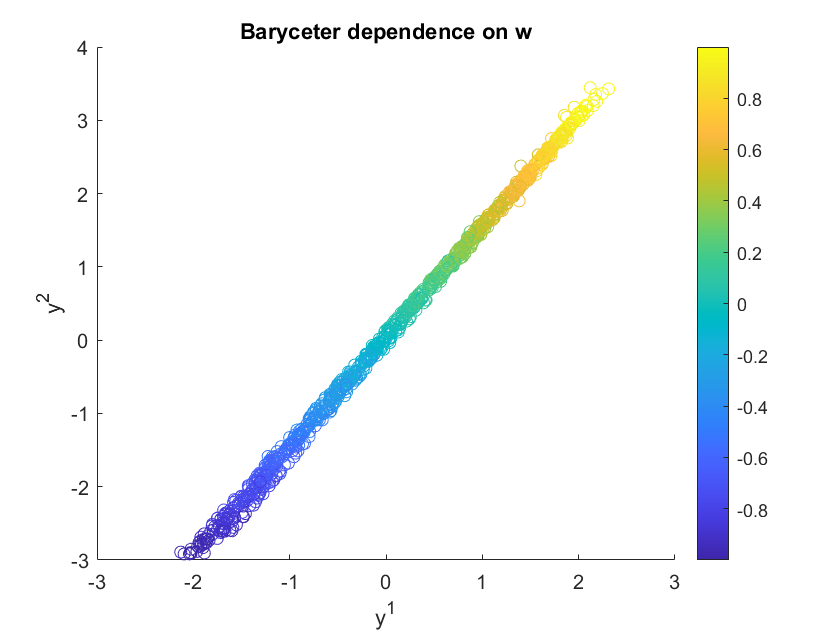}
       \caption{Two-dimensional $y$ dependence on a one-dimensional $w$ (denoted with a colorbar)}
       \label{fig::HS2}
\end{figure}

When the hidden signal $w$ is lower dimensional than $x$, it follows that $y$ must lie in a lower dimensional manifold of $\mathcal{X}$.
Consider an example where $z\in \R^{2}$ with $z\sim\mathcal{N}(0,I)$, $w\in \R$ with $w\sim \mathcal{N}(0, 1)$ and $x = [8z^1z^2+2w, 2z^1+8z^2+3w]$. As shown in Figure \ref{fig::HS2}, after solving the barycenter problem, the resulting $y$ lies on a 1-D manifold that is parametrized (and therefore completely explained) by $w$, which is not generally the case for residual analysis.

\subsection{Hidden patterns in ground-level atmospheric temperature}\label{sec:groundatmotemp}

We consider next the hourly ground-level temperature in Ithaca, NY from 2007 to 2023. The data, available from \href{https://www1.ncdc.noaa.gov/pub/data/uscrn/products/hourly02/?C=D;O=D}{National Oceanic and Atmospheric Administration}, is displayed in the top panel of Figure \ref{fig::time_series_analysis}. We will use the OTBP to investigate the dependence of this temperature on the diurnal and seasonal cycles, and to uncover hidden signals at the synoptic weather and multi-year scales.

\begin{figure}[h!]
\centering
\subfloat[Using only day of year as a factor ]{
\includegraphics[width=0.5\textwidth]{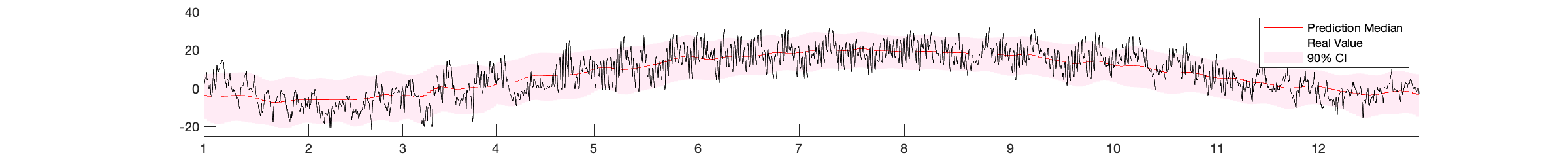}}
\subfloat[Using only hour of day as a factor]{
\includegraphics[width=0.5\textwidth]{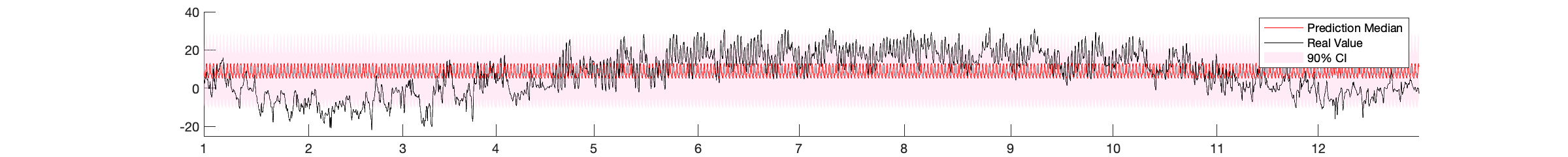}}\\
\subfloat[Using the two periodic time factors together ]{
\includegraphics[width=0.5\textwidth]{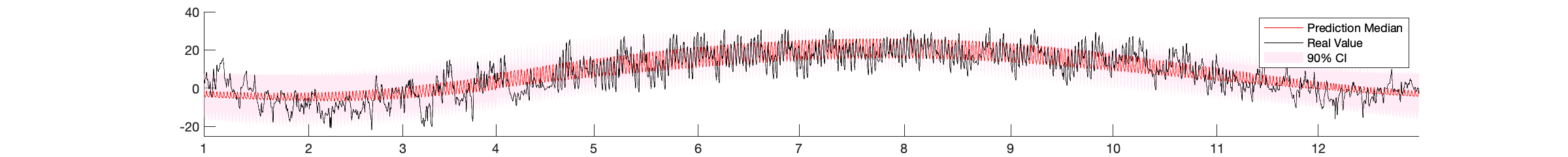}}
\caption{Median temperature and 90\% confidence interval as a function of day of year, time of day and both, displayed over the true temperature for 2007.}
\label{fig::compare_three}
\end{figure}

We first solve the OTBP for $\rho(x | z_1)$, where $x$ is the hourly temperature in Ithaca and $z_1 \in [0, 365.25]$ is the day of the year, a continuous, periodic factor.  Panel (a) of Figure \ref{fig::compare_three} displays the corresponding day-dependent median value of the simulated temperature $X(:, z_1)$ and the corresponding conditional 90\% confidence interval, capturing seasonal effects, superimposed for reference on the true observed temperatures for the year 2007. We use the conditional median and confidence intervals rather than the conditional mean and standard deviation of $x$ because they are more robust statistics and they are also much cheaper to compute: while computing the mean involves averaging $X(y, z)$ over all $\{y_i\}$ for each value of $z$, and similarly for the variance, the monotonicity of $X$ implies that the conditional median of $x$ is just $X(\bar{y}, z)$, where $\bar{y}$ is the median of the $\{y_i\}$), and similarly for confidence intervals, since conditional percentiles of $\rho(x|z)$ translate directly from the corresponding percentiles in $y$. Next we consider instead $\rho(x | z_2)$, where $z_2 \in [0, 24]$ is the time of the day, another continuous and periodic factor, displaying in panel (b) the simulated median diurnal cycle and corresponding confidence interval together with the true $x$ for  2007. Then we combine the two factors and consider $\rho(x | z_1, z_2)$, with results displayed on panel (c) for the the full year 2007. One may notice in all panels how the 90\% confidence interval depends on $z$, often adopting asymmetric shapes around the median and displaying interesting contrasts between day and night. This is one manifestation of the power of capturing the full conditional distribution $\rho(x|z)$, as opposed to just a few statistics, such as the conditional mean value computed in regression. 

\begin{figure}[h!]
\centering
\includegraphics[width=0.2\textwidth]{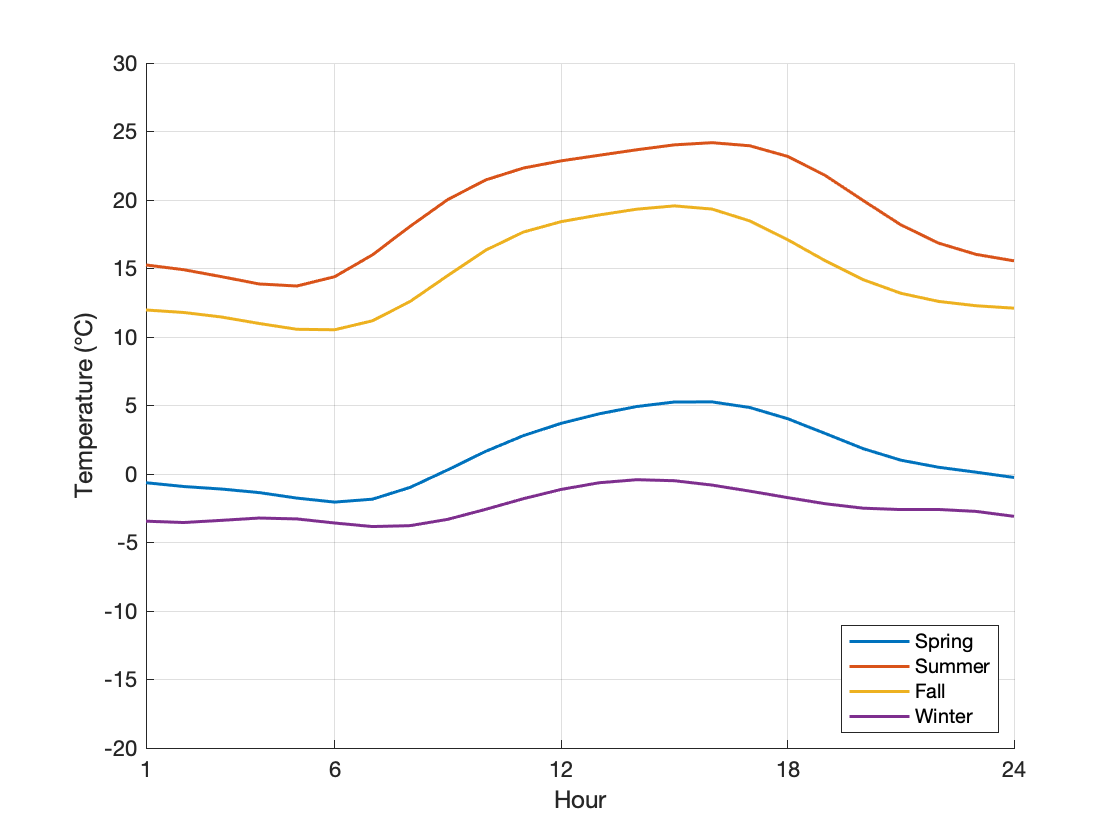}
\hspace{0.02\textwidth}
 \includegraphics[width=0.2\textwidth]{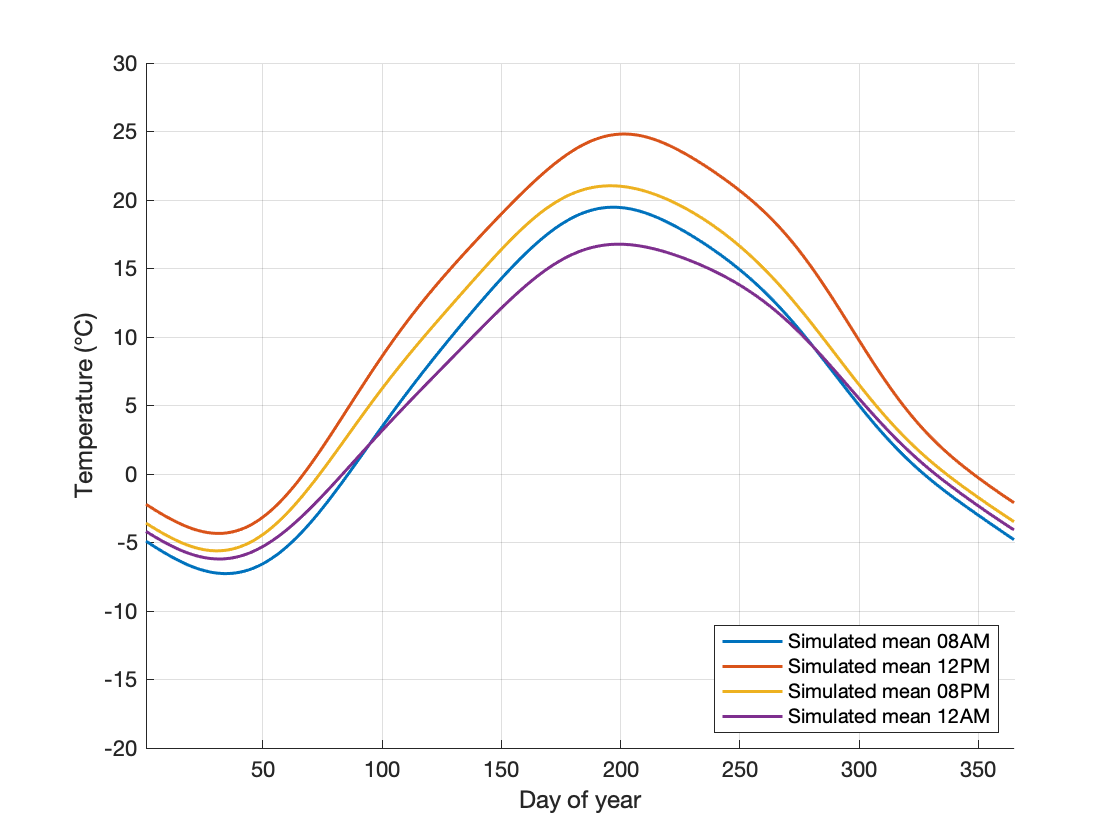} 
 \caption{Median diurnal cycle and seasonal cycle of the temperature in Ithaca, NY, displayed for four days of the year (one per season) and four hours of the day respectively.}
 \label{fig::trend_discovery}
\end{figure}

\begin{figure}[h!]
\centering
\subfloat{
\includegraphics[width=0.5\textwidth]{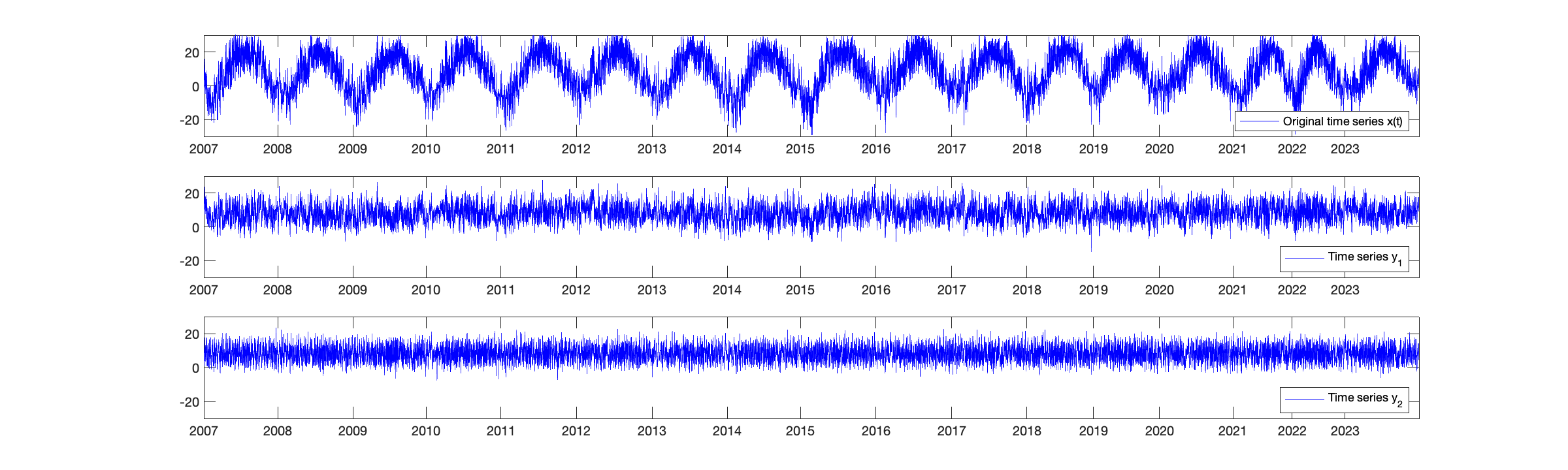}}
\subfloat{
\includegraphics[width=0.5\textwidth]{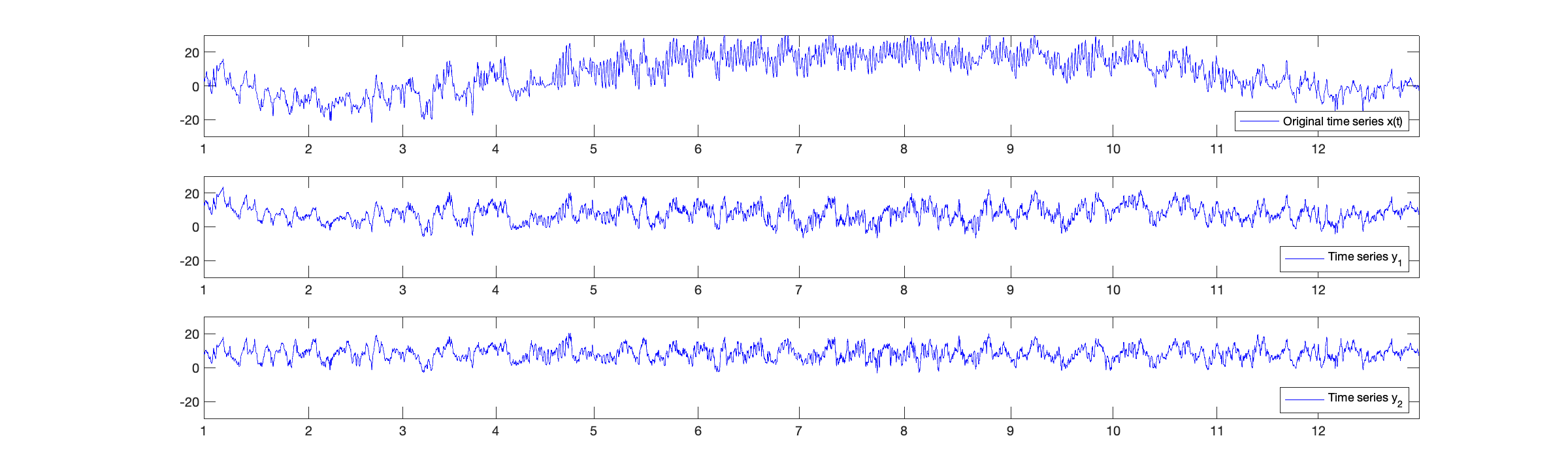}}\\
\caption{The original time series $x(t)$, the $y_1(t)$ resulting from removing the effects of periodic time factors (time of day and day of year), and the $y_2(t)$ resulting from further removing from $y_1$ the synoptic weather and multi-year signals. 
On the right panels, a zoom-in version restricted to 2007.}
\label{fig::time_series_analysis}
\end{figure}

We can see how the diurnal cycle changes over the year not only in mean but also in amplitude and shape. This is seen more clearly in the left panel of Figure \ref{fig::trend_discovery}, displaying the median diurnal cycle for four specific days of the year, corresponding to the  solstices and equinoxes. We can see in detail, for instance, how the Winter Solstice  day is colder, shorter and has smaller day/night contrast than its summer counterpart, and how the day at the Spring Equinox, despite having exactly the same duration as the one at the Fall Equinox, is much colder, has smaller amplitude and a slightly different shape. The right panel of Figure \ref{fig::trend_discovery} similarly shows how the median seasonal cycle depends on the time of the day at which it is considered.

\begin{figure}[h!]
\centering
\subfloat[Median multi-year trend]{
\includegraphics[width=0.5\textwidth]{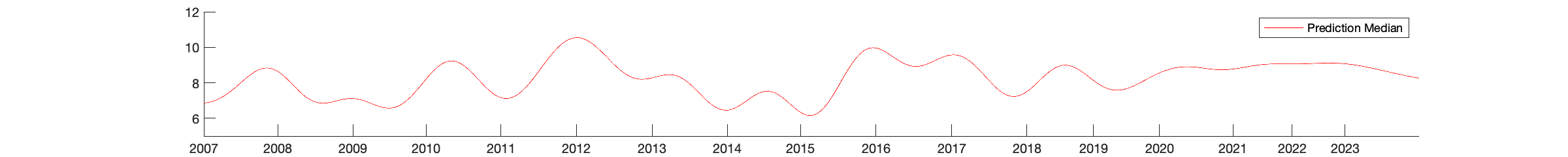}}
\subfloat[Median synoptic weather trend (2007)]{
\includegraphics[width=0.5\textwidth]{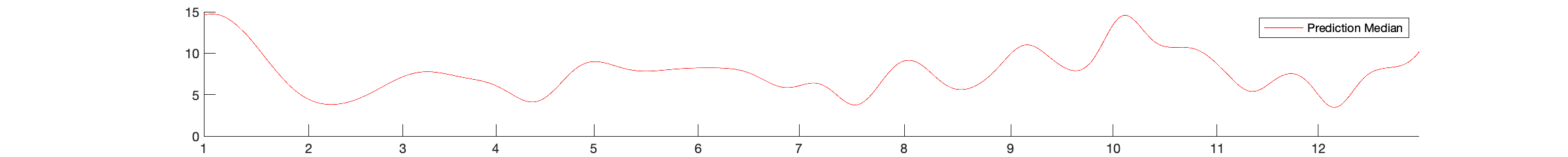}}\\
\subfloat[Median temperature and 90\% confidence interval predicted by $z_3$ and $z_4$ together]{
\includegraphics[width=0.5\textwidth]{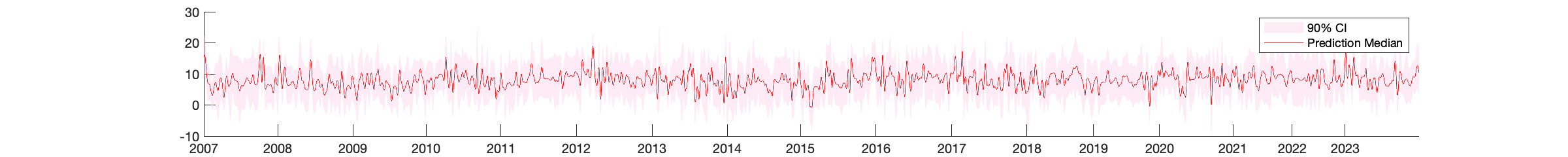}}
\subfloat[Zoom for 2007 superposed on $y_1(t)$]{
\includegraphics[width=0.5\textwidth]{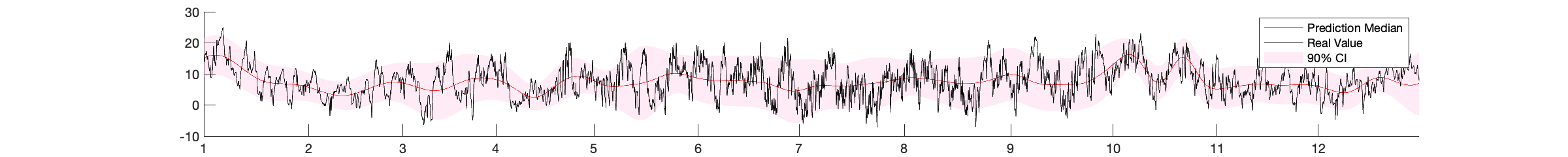}}\\
\subfloat[Median temperature and 90\% confidence interval as a functions of day of year, time of day, multi-year trend and the synoptic weather trend, displayed over the true temperature for 2007]{
\includegraphics[width=0.5\textwidth]{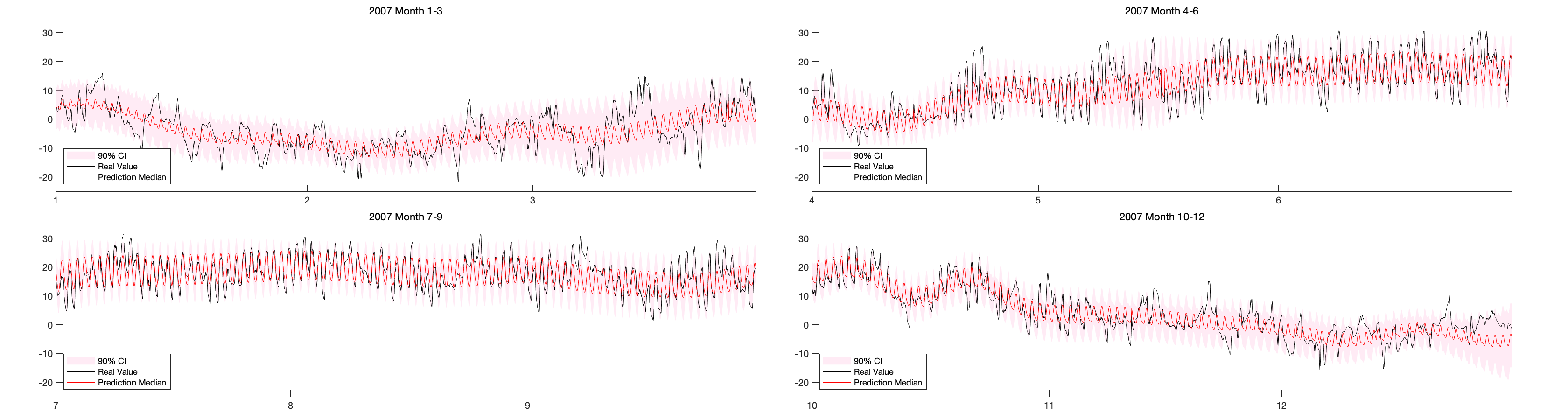}}\\
\caption{Reconstruction of the conditional median temperature $\bar{x}(z)$ and 90\% confidence interval as functions of different combinations of the factors $z_{1,2,3,4}$.}
\label{fig::time_series_analysis_2}
\end{figure}

We switch next to consider the variability of $x$ not explained by $(z_1, z_2)$, as captured in $y$. In order to analyze both synoptic weather and short-term multi-year variability, we introduce two new time factors, $z_3$ and $z_4$, built by rescaling time using different scales. In other words, both $z_3$ and $z_4$ consist just of the time $t$ (measured in hours), but the bandwidths used for the corresponding kernels are of the order of 30 and 3 days, respectively (we set these scales through the parameter $\gamma_z$ defined in the appendix.) As we did for $z_{1,2}$, we remove the variability in $y$ attributable to $z_3$ and $z_4$ separately, then together. Figures \ref{fig::time_series_analysis} and \ref{fig::time_series_analysis_2} show the results of introducing these new factors. Panels (a) and (b) of figure \ref{fig::time_series_analysis} display the original time series $x(t)$, the signal $y_1(t)$ resulting from removing $z_{1,2}$ from $x$ and the signal $y_2(t)$ resulting from removing $z_{3}$ and $z_4$ from $y_1$.
The removal of variability is reflected in the signals' decreasing variance, from $113.71$ for $x$, through $24.71$ for $y_1$, to $17.10$ for $y_2$. Beyond the decreasing variance, one can observe the further explanation of variability in the fact that $y_2$ is much more homogeneous in time than $y_1$, which has a clear inhomogeneity associated with the synoptic weather signal (The fact that similarly $y_1$ does not display the time dependence on the seasonal cycle present on $x$ is far more obvious to the eye.)
Panel (a) of figure \ref{fig::time_series_analysis_2} displays the median temperature dependence on $z_3$, a multi-year signal,  and panel (b) the dependence on $z_4$, corresponding to synoptic weather, over the year 2007. The median temperature and 90\% confidence interval determined by $z_3$ and $z_4$ together is displayed in panel (c) and zoomed-in over 2007 in panel (d). A climate scientist looking at these reconstructions may not only confirm that the method has captured the right scales (a roughly 2-4 year scale for the multiyear signal and around 15 days for the synoptic weather) but also detect individual signals, such as in panel (a) a signal resulting from the El Niño years 2007, 2010 and 2016, and in panel (b) a signal from the North American heat wave of 2007, which may have contributed to the elevated temperature of Ithaca during the late summer and early autumn. Finally, panel (e) displays the full $z_{1,2,3,4}$-dependent  conditional median $\bar{x}(z(t))$ and 90\% confidence interval over 2007. We can see in this plot not only how well the reconstruction has captured the dependence of temperature on time, season and the synoptic and multi-year time scales, but also how it has not captured (by construction) weather signals shorter than a week long. 
Notice that these shorter scales are nonetheless represented as noise in the 90 percentile, a general property of the OTBP methodology: as new factors $z^l$ are introduced, these explain away part of the variability previously present in the conditional distribution $\rho(x|z)$.

 
\subsection{Forecasting of global ocean states}\label{sec:oceanstates}

\begin{figure}[h!]
\centering
\includegraphics[width=0.8\textwidth]{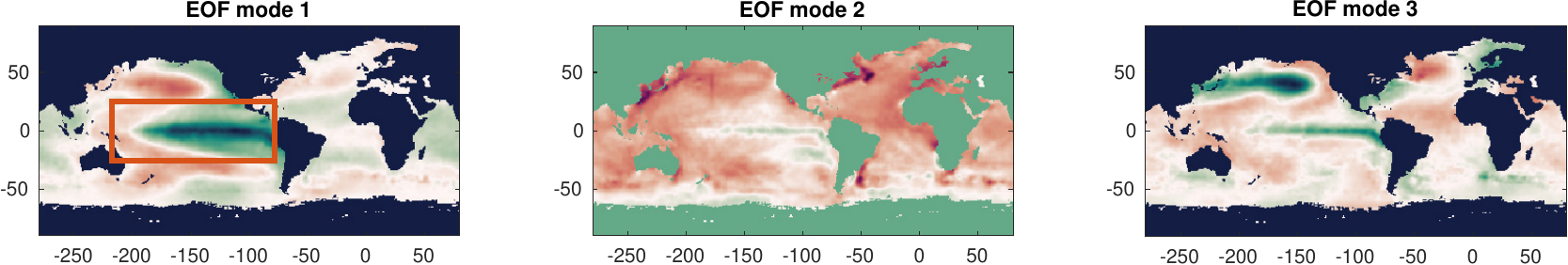}\\
\includegraphics[width=0.8\textwidth]{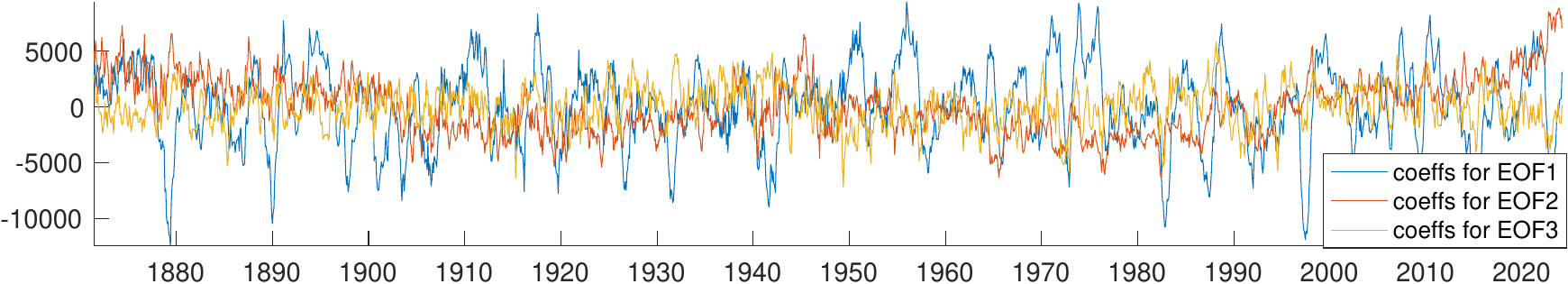}
\caption{Upper panels: $EOF^{1, 2, 3}$.  The red box in $EOF^{1}$ (left most panel) indicates the region of El Ni\~no events. Lower panel: $C_{l}^{1, 2, 3}$ as a function of time $t_{l}$.}
\label{fig::eof}
\end{figure}

We further illustrate the methodology, using it to forecast six months ahead the global sea surface temperature (SST).  The data (available at \href{https://www.metoffice.gov.uk/hadobs/hadisst/data/download.html}{Met Office Hadley Centre observations}) consists of monthly values of the SST from 1870 to 2024, over a global $1\times 1$ latitude-longitude grid.  The resulting dataset $T_{i,j}^{l}$ has a dimension of $180\times 360\times1860$,  corresponding to latitude (indexed by $i$), longitude ($j$) and time ($l$ in months).  The goal is to use historical observations to predict the global SST 6 months into the future. 
In order to extract a lower dimensional time signal from the data,  we apply a standard pre-processing to the whole dataset:
\begin{itemize}
\item Filter out spatial grid points that either lie over land (where SST is undefined) or contain missing data, resulting in $31{,}094$ valid spatial grid points; all subsequent analysis, including EOF computation, is restricted to this filtered spatial domain;
\item De-trend by fitting a linear function of the temporal variable;
\item Explain away the seasonal cycle by removing the mean value at each day of the year from each point on the spatial grid, reducing $T_{i,j}^l$ to the anomaly signal $A_{i,j}^{l}=A(x_{i,j},t_{l})$;
\item Obtain through principal component analysis the first $K$ empirical orthogonal functions of the data \cite{Hannachi2007}, 
$$
A_{i,j}^{l} \approx \sum_{k=1}^{K} \sigma_k \text{C}^{k}_{l} \text{ EOF}^{k}_{i,j} ,
$$
where $\text{EOF}^{k}_{i,j}=\text{EOF}^{k}(x_{ij})$ are static, geographically dependent components of the SST profiles,  and the $\text{C}^{k}_{l}=\text{C}^{k}(t_{l})$ capture their magnitude at time $t_{l}$.  The components are sorted according to $\sigma_k$, proportional to the fraction of the variance of $A$ that they explain.
\end{itemize}
(We could replace the standard pre-processing by a far more informative one based on the OTBP methodology itself.) 
It is known \cite{Takahashi2011} that the first EOF component correlates strongly with ENSO events. Figure \ref{fig::eof} depicts  the first 3 $\text{EOF}^{k}(x)$ as well as their temporal coefficients $\text{C}^{1,2,3}(t)$.  The prediction task then reduces to forecasting the coefficients $\text{C}^{k}$ from their lagged observations  $\text{C}^{l}(t-\Delta t)$,  with $\Delta t=6$ months. 

After approximating the original time series by its first $K=50$ EOFs, we split the data into in-sample (145 years,  1870-2014) and out-of-sample (10 years,  2015-2024) sets.  We apply the procedure to each component $\text{C}^k$ independently,  so $x$ is a one dimensional outcome.  The covariate space $z$ is multidimensional,  consisting of two types of factors:  (1) the $\text{C}^l(t - \Delta t)$ with lagged correlation with $\text{C}^k$ of absolute value greater than 0.1, and (2) time-lagged observations of the same component $\text{C}^k$,  with lags of 6, 12, 24, and 36 months.

\begin{figure}[h!]
\centering
\subfloat[The first 3 coefficients $C^{k}(t)$  (black) and our prediction (red for prediction mean, and pink for one standard deviation away from prediction mean). ]{
\includegraphics[width=0.25\textwidth]{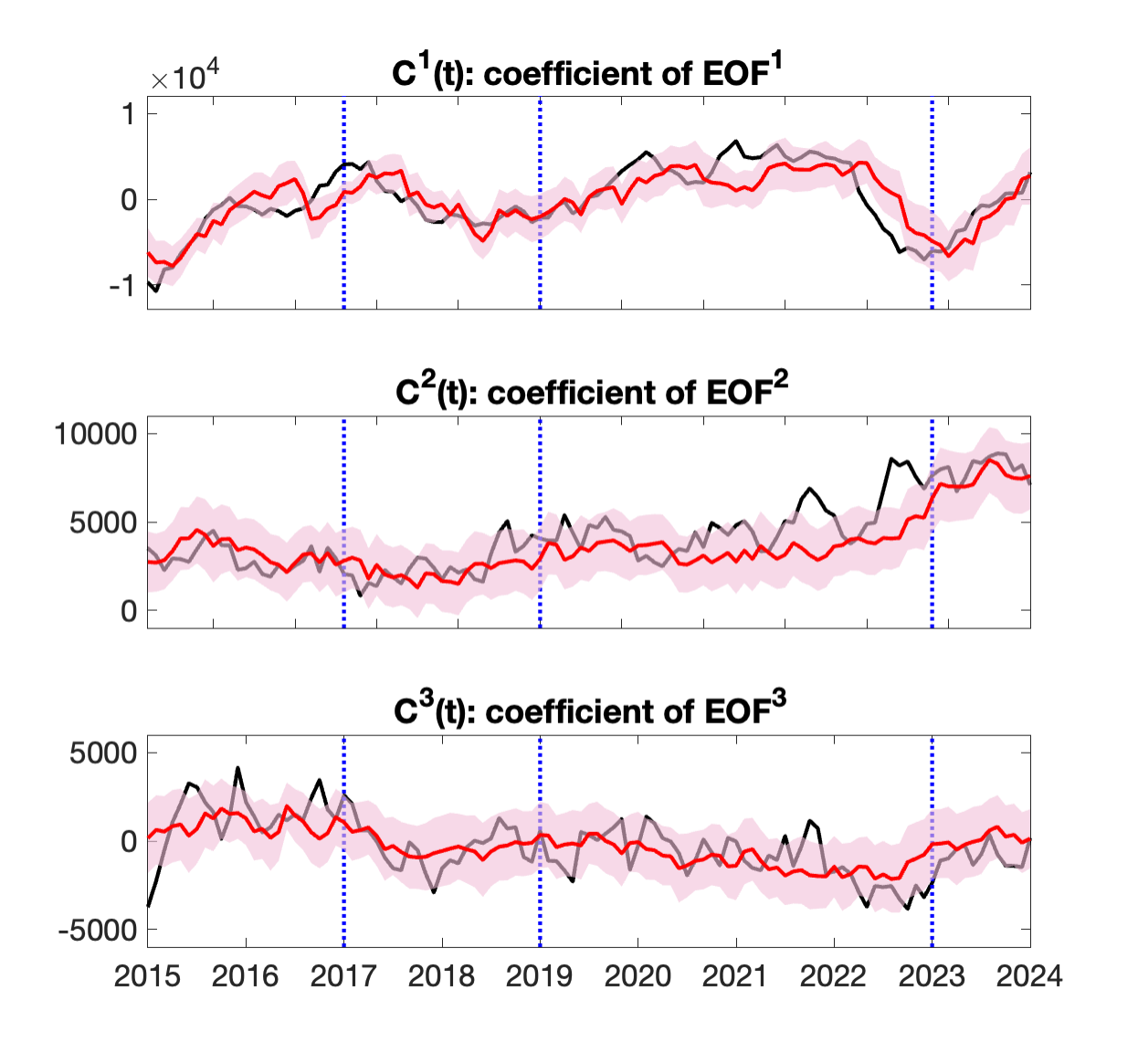}}\hfill
\subfloat[Anomaly at 3 dates (from left to right: strong El Ni\~no,  neutral,  strong La Ni\~na),  color limit: $\pm 1.5 ^\circ C$.  The top row displays the ground truth, while the bottom row contains our prediction, truncated to the top 50 components.  ]{
\includegraphics[width=0.7\textwidth]{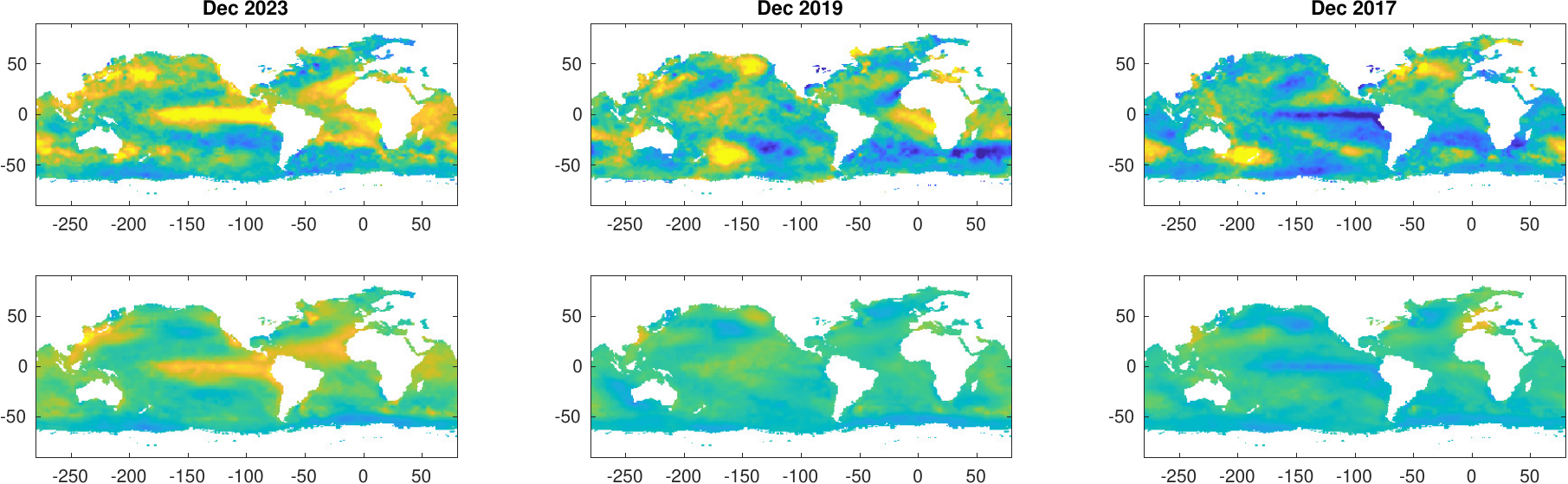}}\\
\caption{Forecast 6 months ahead of the global Sea Surface Temperature: (a) visualization of the signal's first 3 time components, (b) global anomaly.  }
\label{fig::lag6}
\end{figure}

We restrict the family of functions $\mathcal{G}$ and $\mathcal{F}$ in Section \ref{sec:Impl} to include linear and quadratic terms in $y$ and kernels in $z$ space respectively. We cross-validate over the optimal $z$-space kernel bandwidth parameter $\gamma_z$, defined in the appendix.  For each  value of $\gamma_z$ among 40 points uniformly distributed in $[0.2,20]$,  we find the barycenter and compute the $\ell_2$ norm of the difference between true and predicted mean,  evaluated over the out-of-sample data.  The results shown in Figure~\ref{fig::lag6} correspond to the optimal bandwidth that minimizes the norm of difference.

Figure \ref{fig::lag6} depicts the prediction of $A(x,t)$, focusing on SST anomalies for December of 2023,  2019,  and 2017,  corresponding to El Ni\~no,  neutral,  and La Ni\~na events.  %
The predictions performance can be visually inspected in two different spaces: through the prediction of each EOF component $C^k$ (panel (a)) and of $A(x, t)$ for specific times $t$ (panel (b)).
For the second option,  we truncated the prediction to the first 50 components, which explain above 83\% of the variability of the original SST anomaly.  In both spaces,  our method always recovers consistent anomalies both globally and locally within the El Ni\~no region.

\section{Summary and discussion}\label{sec:Summary}

This article develops an efficient methodology for solving the sample-based Monge optimal transport barycenter problem, which takes as input $n$ observed sample pairs $\{x_i, z_i\}$ drawn from an unknown joint distribution $\pi(x, z)$, and produces as output $n$ associated samples $\{y_i\}$ from the barycenter $\mu$ of the conditional distributions $\rho(x|z)$ under $\gamma(z) = \pi(\mathcal{X}, z)$. In addition, it produces $n$ samples $\{x_i^*\}$ drawn from the estimated $\rho(:|z_*)$ for any proposed target value $z_*$ of the covariates $z$ and it estimates $\rho(x|z)$, instrumental for model-free Bayesian inference. A corollary extends the procedure to solve the regular [Monge] optimal transport problem.

Central to the methodology and its applications is a formulation of the OTBP not in terms of the barycenter $\mu$ itself but of the underlying random variable $y = T(x, z)$, which must be statistically independent of the factors $z$. 
A test-based formulation of independence through the uncorrelation between all functions $\{f(z), g(y)\}$ within suitable functional spaces $\{\mathcal{F}, \mathcal{G}\}$ provides an adversarial formulation of the OTBP. Since the best adversarial functions $f$ and $g$ can be found exactly in terms of the first principal components of a matrix $A(y)$, the problem reduces to a single minimization over the map $T$. Solving this problem through gradient descent over $y$ yields a flow that transports each $y_i$ from $x_i$ to $T(x_i, z_i)$. The resulting map $T$ can be inverted in closed form, which facilitates both the simulation and the estimation of $\rho(x|z)$. A byproduct of this closed-form inversion is the extraction of factors $\{f^k\}(z)$ that encode the dependence of $x$ on $z$.

Numerical examples illustrate the applicability of the Monge OTBP and the effectiveness of the methodology proposed. These examples range from synthetic demonstrations of density estimation and simulation, model-free Bayesian inference and hidden signal discovery, to real data applications to weather and climate. Within this article, the latter two are intended only as illustrations of the algorithm at work. A more in-depth study, which requires further extensions of the methodology, should be pursued in field-specific contexts.

This article lays the methodology's general framework. Much more can be done regarding the adaptive choice of the functional spaces $\mathcal{F}$ and $\mathcal{G}$, for which we have proposed here just a handful of simple choices. Since the range of options to explore on adaptive functional spaces is too broad for a single article,  exploring them further here would take us too far afield. We also choose not to dwell in this article on other extensions, such as going beyond gradient descent, as required for factor discovery, or building functional spaces $\mathcal{G}$ better-suited for high-dimensional outcome spaces $\mathcal{X}$. We believe that the proposed methodology can be extended in a number of meaningful directions, making it an effective, robust, versatile and conceptually sound approach to a broad set of tasks in data analysis.

\appendix

\section{Appendix: a data-adapted functional space}

This appendix describes the choice of functional spaces $\mathcal{F}$ and $\mathcal{G}$ used in our numerical examples. Since the two constructions are entirely similar, we describe only the space $\mathcal{F}$. Exploring other, potentially much richer choices of adaptive functional spaces goes beyond this article's scope.

Since the components of $z$ can be of arbitrary type, including real, periodic, categorical and more (not so those of $y$, which are typically real), we first embed $z$ in an Euclidean space $R^k$ as follows. For each component $z^l \not\in R$ of $z$: 
1) when $z^l$ is periodic with period $T$, we embed it in $R^2$, mapping $z^l$ to $w$ on the unit circle,
$ w\left(z^l\right) = \left[\cos\left(\frac{2\pi}{T} z^l\right), \sin\left(\frac{2\pi}{T} z^l\right)\right]$;
2) when $z^l$ is categorical with $h$ discrete values $v_k$, we embed it in $R^{h-1}$, mapping the $\{v_k\}$ to $h$ equidistant points, the vertices of a regular simplex;
3) for variables $z$ of a more complex type, such as images, distributions or graphs, we introduce an application-specific distance among them and embed them into some $R^h$ accordingly.
Having done this, we can restrict attention to $\mathcal{Z} = R^{d_z}$.

For any function $f(z)$ and any probability density $\gamma(z)$, we have
$$\gamma(z) = \int \gamma\left(z'\right) \delta\left(z - z'\right) dz', \quad
f(z) = \int f\left(z'\right) \delta\left(z - z'\right) dz' = 
\int f\left(z'\right) \frac{\delta\left(z - z'\right)}{\int \gamma\left(z''\right) \delta\left(z' - z''\right) dz''}\ \gamma\left(z'\right) dz' $$
(Notice that the last expression is only valid for values of $z$ within the support of the distribution $\gamma$.)  
Mollifying  $\delta(x-y)$ to a non-negative function $K(x, y)$ that concentrates near $x=y$ and integrates to $1$ over $x$ yields
$f(z) \approx 
\int f\left(z'\right) \frac{K\left(z, z'\right)}{\int \gamma\left(z''\right) K\left(z', z''\right) dz''} \gamma\left(z'\right) dz'$,
which in terms of samples $ z^c_j \sim \gamma(z)$ yields the empirical version
\begin{equation}
f(z) \approx 
\sum_j f\left(z^c_j\right) \frac{K\left(z, z^c_j\right)}{\sum_h  K\left(z^c_j, z^c_h\right)}. 
\label{KR_not_quite}
\end{equation}
Notice that nothing in our argument requires $K(x, y)$ to be a symmetric function, so it is not a ``kernel'' in the conventional sense. This allows us to use center-dependent bandwidths, coarser where the data is sparser.

If $K$ is a smooth function of $z$, so is the right-hand side of (\ref{KR_not_quite}) for any choice of  $f\left(z^c_j\right)$. We conclude that
$f(z) = 
\sum_j a_j \frac{K\left(z, z^c_j\right)}{\sum_h  K\left(z^c_j, z^c_h\right)}$
parameterizes arbitrary smooth functions of $z$ within the support of $\gamma(z)$.
Since the denominator does not depend on $z$, we can absorb it into the definition of $a_j$, which yields
\begin{equation}
f(z) = 
 F(z) a,  \quad F^j(z) = K\left(z, z^c_j\right),  
 \label{Fa}
 \end{equation}
where it is no longer required that the $\{F^j(z)\}$ integrate to one, since the corresponding normalizing constants can also be absorbed into the $\{a_j\}$.   Then the space $\mathcal{F}$ of smooth functions $f(z)$ agrees with the column space of the operator $F$.
To consider functions with zero mean, it is enough to subtract the mean of each column of $F$.

There is no need for the set of centers $\{z^c_j\}$ for $K$ and the set of points $\{z_i\}$ where $f$ is to be evaluated to agree; it is enough that the support of the distribution $\gamma$ underlying the former contains the support of the latter. 
Using a number $m \ll n$ of centers reduces the computational cost associated to evaluating the kernels. Moreover, when applied to $g(y)$, one needs to decouple the centers $\{y^c_j\}$ from the samples $\{y_i\}$, as only the latter are arguments over which the objective function $L$ is minimized.  
A simple procedure for selecting $m \ll n$ centers is through k-means applied to the $\{z_i\}$, which has been the choice adopted for all examples in Section \ref{sec:NumEx}, where we have set $m = \min([\sqrt{n}], m_{max})$, with $m_{max}$ set by the user. 
We adopted for $K(z, z^c)$ a Gaussian function with center-dependent inverse covariance matrix $S_j = S({z^c_j})$:
$ K\left(z, z^c_j\right) = e^{-\frac{1}{2} \left(z - z_j^c\right)' S_{j} \left(z - z_j^c\right)}$.

Tuning the $\{S_j\}$ is critical for extracting as much dependence of $x$ on $z$ as possible, by capturing the right functions $f(z)$ and $g(y)$. 
The appropriateness of a functional space for $f(z)$ depends not just on the samples $\{z_i\}$ but also on the $\{x_i\}$. 
Yet we do not know the form of this dependence before hand, since determining it is precisely our algorithm's goal. Thus we  make a choice based not on the relation between $x$ and $z$ but on the data available for each.  
 The most natural function $f(z)$ to attempt to capture when looking only at the $\{z_i\}$ is their underlying probability density $\gamma(z)$. We apply the following adaptive procedure to determine the corresponding $\{S_j^0\}$.

Compute first a global empirical mollified covariance matrix $\Sigma$ and its inverse $S_g$,
$$ \Sigma^{kl} = \frac{1}{m} \sum_{i=1}^m \left(\left(z_i^c\right)^k - \bar{z}^k\right) \left(\left(z_i^c\right)^l - \bar{z}^l\right) + \epsilon I, \quad \epsilon = \frac{\hbox{var}(z)}{m}, \quad
\quad S_g = \Sigma^{-1}. $$

Introducing an adjustable parameter $\alpha$, define
$ K_i^j = K_{g}\left(z^c_i, z^c_j\right) = 
e^{-\frac{1}{2 \alpha^2} \left(z^c_i - z^c_j\right)'  S_g \left(z^c_i - z^c_j\right)}$,
estimate $\gamma^i_{\alpha}\left(z^c_i\right)$ through leave-one-out kernel density estimation,
$ \gamma^i_{\alpha}\left(z^c_i\right) \propto  \frac{1}{\alpha^d}\sum_{j \ne i} K_i^j$,
determine $\alpha$ through leave-one-out maximal likelihood,
$ \alpha^* = \arg\max_{\alpha} L = \sum_{i=1}^m \log\left(\gamma^i_{\alpha}\left(z^c_i\right)\right)$,
and define accordingly $S_{\alpha_*} = \frac{S_g}{\alpha_*^2}$.
A practical choice is to maximize $L$ over a finite set of candidate $\alpha$'s centered around the rule-of-thumb value
$ \alpha_{r.o.th.} = \frac{\left(\frac{4}{d+2}\right)^{\frac{1}{d+4}}}{m^{\frac{1}{d+4}}}$.
Rescale the $S_{\alpha_*}$ locally using the estimated $\gamma^i_{\alpha_*}$,
$ S_j = \left(\gamma_{\alpha_*}\left(z^c_j\right)\right)^{\frac{2}{d}} S_{\alpha_*}$,
so that the the number of points $\{z^c_i\}$ within the effective support of the corresponding function $K(z, z^c_j)$ is roughly independent of $j$. 
Finally, introducing a new global adjustable parameter $\beta$, write
$ K_j^{\beta}(z) = e^{-\frac{1}{2 \beta^2} \left(z - z^c_j\right)'  S_j \left(z - z^c_j\right)}$,
with $\beta$ determined again through leave-one-out maximal likelihood, and define
$ S^0_j = \frac{S_j}{\beta_*^2}$:
$$ \beta_* = \arg\max_{\beta} L = \sum_{i=1}^m \log\left(\gamma_{\beta}\left(z^c_i\right)\right), \quad \gamma_{\beta}\left(z^c_i\right) \propto  \frac{1}{\beta^d}\sum_{j \ne i} \frac{K_i^j}{\gamma_{\alpha_*}\left(z^c_j\right)}, \quad K_i^j = K^{\beta}_{j}\left(z^c_i\right), $$

There are at least two reasons why we may want to add one or more free parameters to the determination of the bandwidths. One is that the ideal $S_j$ should depend not just on the $\{z_i\}$ but also on their relation to the $\{x_i\}$, and one straightforward way to address this dependence is through cross-validation over such free parameters. The second reason is that often the bandwidths are determined not by the data alone but also from the scales that one seeks to resolve, as in the multi time-scale analysis of ground temperature of Section \ref{sec:groundatmotemp}. In view of this, we divide $S^0_j$ by an externally provided constant ${\gamma_z}^2$, which we either cross-validate over or set based on the scales that we seek to resolve.

\section*{Acknowledgments}
The work by Y. Wang and E. G. Tabak was partially supported by a grant from the ONR, and W. Zhao's research was supported in part by the Pacific Institute for the Mathematical Sciences and the Simons Foundation.

\end{document}